\documentstyle[a4wide]{article}
\begin{document}
\newcommand{\gna}{
\begin{picture}(40,12)(0,6)
\put(0,10){\vector(1,0){39}}
\put(20,15){\circle*{3}}
\put(20,5){\makebox(0,0){$\Gamma$}}
\end{picture}}

\newcommand{\fna}{
\begin{picture}(40,12)(0,6)
\put(0,10){\vector(1,0){39}}
\put(20,15){\circle*{3}}
\put(20,5){\makebox(0,0){$\Phi$}}
\end{picture}}

\newcommand{\pna}{
\begin{picture}(40,12)(0,6)
\put(0,10){\vector(1,0){39}}
\put(20,15){\circle*{3}}
\put(20,5){\makebox(0,0){$\Psi$}}
\end{picture}}

\newcommand{\fpna}{
\begin{picture}(40,12)(0,6)
\put(0,10){\vector(1,0){39}}
\put(20,15){\circle*{3}}
\put(20,5){\makebox(0,0){$\Phi\cdot\Psi$}}
\end{picture}}

\newcommand{\pfna}{
\begin{picture}(40,12)(0,6)
\put(0,10){\vector(1,0){39}}
\put(20,15){\circle*{3}}
\put(20,5){\makebox(0,0){$\Phi \Psi$}}
\end{picture}}

\newcommand{\qed}{\nopagebreak
\begin{flushright} {\em q.e.d.} \end{flushright}}
\newcommand{\br}{\mbox{$\overrightarrow{\mbox{\bf b}}$}}
\newcommand{\bl}{\mbox{$\overleftarrow{\mbox{\bf b}}$}}
\newcommand{\kr}{\mbox{$\overrightarrow{\mbox{\bf k}}$}}
\newcommand{\kl}{\mbox{$\overleftarrow{\mbox{\bf k}}$}}
\newcommand{\ms}{\mbox{\boldmath{${\sigma}$}}}
\newcommand{\md}{\mbox{\boldmath{${\delta}$}}}
\newcommand{\mal}{\mbox{\boldmath{${\alpha}$}}}
\newcommand{\mbe}{\mbox{\boldmath{${\beta}$}}}
\newcommand{\mep}{\mbox{\boldmath{${\epsilon}$}}}
\newcommand{\met}{\mbox{\boldmath{${\eta}$}}}
\newcommand{\mtv}{\mbox{\boldmath{$T$}}}
\newcommand{\tv}{\mbox{\rm T}}
\newcommand{\mfv}{\mbox{\boldmath{$F$}}}
\newcommand{\mgv}{\mbox{\boldmath{$G$}}}
\newcommand{\mj}{\mbox{\bf 1}}
\newcommand{\ccc}{\mbox{\bf CartCl}}
\newcommand{\frd}{\mbox{\bf Finord}}
\newcommand{\cc}{\mbox{\bf Cart}}
\newcommand{\sym}{\mbox{\bf SyMon}}
\newcommand{\symcl}{\mbox{\bf SyMonCl}}
\newcommand{\skup}{\mbox{\bf Set}}
\newcommand{\mc}{\mbox{\bf c}}
\newcommand{\mw}{\mbox{\bf w}}
\newcommand{\mk}{\mbox{\bf k}}
\newcommand{\mb}{\mbox{\bf b}}
\newcommand{\ri}{\mbox{\rm I}}
\newcommand{\mi}{\mbox{\scriptsize{\rm I}}}
\newcommand{\mpm}{\mbox{\boldmath{$\pi$}}}           
\newcommand{\dkz}{\\[0.3cm] {\bf proof}\hspace{2em}}
\newcommand{\str}{\rightarrow}

\baselineskip=1.05\baselineskip

\title 
{Coherence in Substructural Categories}

\author{
Zoran Petri\' c
\\Faculty of Mining and Geology
\\University of Belgrade
\\Dju\v sina 7
\\11000 Belgrade, Yugoslavia
\\e-mail: zpetric@mi.sanu.ac.yu
 }

\date{ }
\maketitle
\begin{abstract}
\noindent It is proved that MacLane's coherence results
for monoidal and symmetric monoidal categories can be
extended to some other categories with multiplication;
namely, to relevant, affine and cartesian categories. All
results are formulated in terms of natural transformations
equipped with ``graphs'' (g-natural transformations) and
corresponding morphism theorems are given as consequences.
Using these results, some basic relations between the free
categories of these classes are obtained.
\end{abstract}
\vspace{1em}
\noindent In \cite{natass} MacLane has shown that monoidal and
symmetric monoidal categories are coherent, although the complete definition
of the notion was given for the first time in \cite{clcat}.
Strictly keeping to that definition, we show that relevant, affine and 
cartesian categories are coherent. All the categories above we call
{\em substructural} because they correspond to the minimal fragments of
Associative Lambek's Calculus, linear, relevant, BCK and intuitionistic
logic that are sufficient to describe the underlying structural rules
(see \cite{sl} to find about different aspects
of substructural logics). We use equational axiomatizations of these
categories, which  originate from \cite{mfc}, rather than postulating
the commutativity of certain diagrams. Of course, one who prefers
diagram-chasing can easily convert these equations into commutative diagrams.
\section{Substructural categories}

\noindent By {\em category with multiplication}
we mean a category $\cal A$
together with a bifunctor $\cdot : {\cal A} \times {\cal A} \str {\cal A}$
and a special object I.  Categories with multiplication can be
axiomatized by postulating the following {\em equations} between arrows
\begin{tabbing}
{\em categorial equations}
\\
$(cat{\mbox{\bf 1}})$ {\mbox{\hspace{1cm}}}
 \= $f{\mbox{\bf 1}}_{A}=f={\mbox{\bf 1}}_{B}f$ for all $f:A \str B$
\\
$(cat{\mbox{\bf 2}})$ 
\> $h(gf)=(hg)f$ for all $f,g,h \in Mor({\cal A})$
\end{tabbing}
\begin{tabbing}
{\em functorial equations}
\\
$({\cdot})$ {\mbox{\hspace{1.5cm}}}
\= $(g_{1}f_{1}){\cdot}(g_{2}f_{2})=(g_{1}{\cdot}g_{2})(f_{1}{\cdot}f_{2})$
\\
$({\cdot}{\mbox{\bf 1}})$ 
\> ${\mbox{\bf 1}}_{A}{\cdot}{\mbox{\bf 1}}_{B}={\mbox{\bf 1}}_{A{\cdot}B}$
\end{tabbing}
A category with multiplication is {\em monoidal} if there are special arrows
for all objects $A$, $B$ and $C$ 
\[{\mbox{\boldmath $\sigma$}}_{A}:{\mbox{\rm I}}{\cdot}A \str A
 {\mbox{\hspace{1cm}}}
\md_{A}:A{\cdot}{\mbox{\rm I}} \str A\]
\[\ms_{A}^{i}:A \str {\mbox{\rm I}}{\cdot}A {\mbox{\hspace{1cm}}}
\md_{A}^{i}:A \str A{\cdot}{\mbox{\rm I}}\]
\[\br_{A,B,C}:A{\cdot}(B{\cdot}C) \str (A{\cdot}B){\cdot}C
{\mbox{\hspace{1cm}}}
\bl_{A,B,C}:(A{\cdot}B){\cdot}C \str A{\cdot}(B{\cdot}C)\]
and if it satisfies
\begin{tabbing}
\ms\md-{\em equations}
\\
(\ms){\mbox{\hspace{1cm}}}\= For $f:A \str B$,{\mbox{\hspace{1cm}}}
\= $f\ms_{A}=\ms_{B}(\mj_{\mi}{\cdot}f)$.
\\
(\md) \> For $f:A \str B$,
\> $f\md_{A}=\md_{B}(f{\cdot}\mj_{\mi})$.
\\
$(\ms\ms^{i})$
\> $\ms_{A}\ms_{A}^{i}=\mj_{A}$,
\> $\ms_{A}^{i}\ms_{A}=\mj_{\mi{\cdot}A}$
\\
$(\md\md^{i})$
\> $\md_{A}\md_{A}^{i}=\mj_{A}$,
\> $\md_{A}^{i}\md_{A}=\mj_{A{\cdot}\mi}$
\\
$(\ms\md)$
\> $\ms_{\mi}=\md_{\mi}$
\end{tabbing}
\begin{tabbing}
{\bf b}-{\em equations}
\\
({\bf b}){\mbox{\hspace{0.5cm}}}\= For $f:A \str D$, $g:B \str E$ and
$h:C \str F$, \= $((f{\cdot}g){\cdot}h)\br_{A,B,C}=\br_{D,E,F}
(f{\cdot}(g{\cdot}h))$.
\\
({\bf b}{\bf b}) \> $\br_{A,B,C}\bl_{A,B,C}=
\mj_{(A{\cdot}B){\cdot}C}$,
\> $\bl_{A,B,C}\br_{A,B,C}=\mj_{A{\cdot}(B{\cdot}C)}$
\\(\ms\md{\bf b}) \> $(\md_{A}{\cdot}\mj_{B})
\br_{A,\mi,B}=\mj_{A}{\cdot}\ms_{B}$
\\
({\bf b}5)
\> $\br_{A{\cdot}B,C,D}\br_{A,B,C{\cdot}D}=
(\br_{A,B,C}{\cdot}\mj_{D})\br_{A,B{\cdot}C,D}(\mj_{A}
{\cdot}\br_{B,C,D})$
\end{tabbing}

\noindent A monoidal category is {\em symmetric monoidal}
if it has the special arrow
\[ \mc_{A,B}:A{\cdot}B \str B{\cdot}A \]
for every pair $(A,B)$ of its objects,
and if the following equations hold
\begin{tabbing}
\mc-{\em equations}
\\
(\mc){\mbox{\hspace{1cm}}} \= For $f:A \str C$ and $g:B \str D$,
$(g{\cdot}f)\mc_{A,B}=\mc_{C,D}(f{\cdot}g)$
\\
(\mc\mc) \> $\mc_{B,A}\mc_{A,B}=\mj_{A{\cdot}B}$
\\
(\ms\md\mc) \> $\ms_{A}\mc_{A,\mi}=
\md_{A}$
\\
({\bf b}\mc6)
\> $\br_{C,A,B}\mc_{A{\cdot}B,C}\br_{A,B,C}=
(\mc_{A,C}{\cdot}\mj_{B})\br_{A,C,B}(\mj_{A}{\cdot}\mc_{B,C})$
\end{tabbing}

\noindent A symmetric monoidal category is {\em relevant} if it has
the special arrow
\[ \mw_{A}:A \str A{\cdot}A \]
for every object $A$, and if the following equations hold
\begin{tabbing}
\mw-{\em equations}
\\
(\mw){\mbox{\hspace{1cm}}} \= For $f:A \str B$,
$(f{\cdot}f)\mw_{A}=\mw_{B}f$.
\\
(\ms\md\mw)
\> $\ms_{\mi}\mw_{\mi}=\mj_{\mi}$
\\
({\bf b}\mw)
\> $\br_{A,A,A}(\mj_{A}{\cdot}\mw_{A})\mw_{A}=(\mw_{A}{\cdot}
\mj_{A})\mw_{A}$
\\
(\mc\mw) \> $\mc_{A,A}\mw_{A}=\mw_{A}$
\\
({\bf b}\mc\mw8)
\>$\mc^{m}_{A,B,A,B}\mw_{A{\cdot}B}=\mw_{A}{\cdot}\mw_{B}$,
where
\\
\> $\mc^{m}_{A,B,C,D}=^{df}\br_{A,C,B{\cdot}D}
(\mj_{A}{\cdot}(\bl_{C,B,D}(\mc_{B,C}{\cdot}\mj_{D})
\br_{B,C,D}))\bl_{A,B,C{\cdot}D}$ 
\end{tabbing}

\noindent A symmetric monoidal category is {\em affine} if it has
the special arrow
\[ \mk_{A}:A \str {\mbox{\rm I}} \]
for every object $A$ and if the following equations hold
\begin{tabbing}
\mk-{\em equations}
\\
(\mk){\mbox{\hspace{1cm}}} \= For $f:A \str B$, $\mk_{A}=\mk_{B}f$
\\
(\mj\mk) \> $\mk_{\mi}=\mj_{\mi}$
\end{tabbing}

\noindent If a symmetric monoidal category
is both relevant and affine (described
in the same language) and if its arrows satisfy the following equations
\[ (\ms\mk\mw){\mbox{\hspace{0.5cm}}} \ms_{A}(\mk_{A}{\cdot}
\mj_{A})\mw_{A}=\mj_{A},{\mbox{\hspace{1cm}}}
(\md\mk\mw){\mbox{\hspace{0.5cm}}} \md_{A}(\mj_{A}{\cdot}\mk_{A})
\mw_{A}=\mj_{A} \]
then we say it is {\em cartesian}.
\\
We call this axiomatization of cartesian categories  
{\em structural-equational}. It differs from the standard
equational axiomatization (see \cite{ithocl}) of these categories.
The latter is based on the universality of product and uses as primitive, arrows
$\mj_{A}:A \str A$,
$\mpm_{A,B}:A{\cdot}B \str A$, $\mpm_{A,B}':A{\cdot}B \str B$ and
$\mk_{A}:A \str \ri$ for all objects $A$ and $B$, and a partial
binary operation on arrows $\langle \; , \; \rangle$, such that
  \[ \frac{f:C \str A {\mbox{\hspace{1cm}}} g:C \str B}
  {{\langle}f,g{\rangle}:C \str A{\cdot}B}\]
Equations that hold are the categorial equations plus
\begin{tabbing}
({\rm E2}){\mbox{\hspace{1cm}}} \= $f=\mk_{A}$, for every $f:A \str \ri$.
\\
({\rm E3a.}) \> $\mpm_{A,B}{\langle}f,g{\rangle}=f$,
for $f:C \str A$ and $g:C \str B$.
\\
({\rm E3b.}) \> $\mpm_{A,B}'{\langle}f,g{\rangle}=g$,
for $f:C \str A$ and
$g:C \str B$.
\\
({\rm E3c.}) \> ${\langle}\mpm_{A,B}h,\mpm_{A,B}'h{\rangle}
=h$, for $h:C \str A{\cdot}B$.
\end{tabbing}

To show that these two axiomatizations are extensionally equivalent we have
to define 
\[ \mpm_{A,B}=_{df}\md_{A}(\mj_{A}{\cdot}\mk_{B}),{\mbox{\hspace{0.5cm}}}
\mpm_{A,B}'=_{df}\ms_{B}(\mk_{A}{\cdot}\mj_{B}), \]
and for $f:C \str A$ and $g:C \str B$
\[{\langle}f,g{\rangle}=_{df}(f{\cdot}g)\mw_{C},\]
in structural case. Then it is easy to show that {\rm (E2)--(E3c.)} hold.

Conversely, if we start with the standard axiomatization,
then we can define
\[ \ms_{A}=_{df}\mpm_{\mi,A}',{\mbox{\hspace{0.5cm}}}
\ms_{A}^{i}=_{df}{\langle}\mk_{A},\mj_{A}{\rangle}, \]
\[ \md_{A}=_{df}\mpm_{A,\mi},{\mbox{\hspace{0.5cm}}}
\md_{A}^{i}=_{df}{\langle}\mj_{A},\mk_{A}{\rangle}, \]
\[ \br_{A,B,C}=_{df}{\langle}{\langle}\mpm_{A,B{\cdot}C},
\mpm_{B,C}\mpm_{A,B{\cdot}C}'{\rangle},\mpm_{B,C}'\mpm_{A,B{\cdot}C}'
{\rangle} \]
\[ \bl_{A,B,C}=_{df}{\langle}\mpm_{A,B}\mpm_{A{\cdot}B,C},
{\langle}\mpm_{A,B}'\mpm_{A{\cdot}B,C},\mpm_{A{\cdot}B,C}'
{\rangle}{\rangle}, \]
\[ \mc_{A,B}=_{df}{\langle}\mpm_{A,B}',\mpm_{A,B}{\rangle},
{\mbox{\hspace{0.5cm}}}\mw_{A}=_{df}{\langle}\mj_{A},\mj_{A}{\rangle}, \]
and for $f:A \str C$ and $g:B \str D$
\[ f{\cdot}g=_{df}{\langle}f\mpm_{A,B},g\mpm_{A,B}'{\rangle}. \]
It is straightforward to prove that functorial, \ms\md, {\bf b},
\mc, \mw, \mk- equations, as well as (\ms\mk\mw) and (\md\mk\mw) hold.

In addition we have to prove that the ``double translation'' will take us
to the same notions, i.e., to show that in the structural
axiomatization the following equations hold
\[ \ms_{A}=\ms_{A}(\mk_{\mi}{\cdot}\mj_{A}),{\mbox{\hspace{0.5cm}}}
\ms_{A}^{i}=(\mk_{A}{\cdot}\mj_{A})\mw_{A} \]
\[ \md_{A}=\md_{A}(\mj_{A}{\cdot}\mk_{\mi}),{\mbox{\hspace{0.5cm}}}
\md_{A}^{i}=(\mj_{A}{\cdot}\mk_{A})\mw_{A} \]
\begin{eqnarray*}
\br_{A,B,C}&=&(((\md_{A}(\mj_{A}{\cdot}\mk_{B{\cdot}C})){\cdot}
(\md_{B}(\mj_{B}{\cdot}\mk_{C})\ms_{B{\cdot}C}(\mk_{A}{\cdot}
\mj_{B{\cdot}C})))\mw_{A{\cdot}(B{\cdot}C)}){\cdot}\\
& & (\ms_{C}(\mk_{B}{\cdot}\mj_{C})\ms_{B{\cdot}C}
(\mk_{A}{\cdot}\mj_{B{\cdot}C})))\mw_{A{\cdot}(B{\cdot}C)}\\
\bl_{A,B,C}&=&((\md_{A}(\mj_{A}{\cdot}\mk_{B})
\md_{A{\cdot}B}(\mj_{A{\cdot}B}{\cdot}\mk_{C})){\cdot}
(((\ms_{B}(\mk_{A}{\cdot}\mj_{B})\md_{A{\cdot}B}
(\mj_{A{\cdot}B}{\cdot}\mk_{C})){\cdot} \\
& & (\ms_{C}(\mk_{A{\cdot}B}{\cdot}\mj_{C})))\mw_{(A{\cdot}B){\cdot}C}))
\mw_{(A{\cdot}B){\cdot}C}\\
\mc_{A,B}&=&((\ms_{B}(\mk_{A}{\cdot}\mj_{B})){\cdot}
(\md_{A}(\mj_{A}{\cdot}\mk_{B})))\mw_{A{\cdot}B}\\
\mw_{A}&=&(\mj_{A}{\cdot}\mj_{A})\mw_{A}\\
f{\cdot}g&=&((f\md_{A}(\mj_{A}{\cdot}\mk_{B})){\cdot}
(g\ms_{B}(\mk_{A}{\cdot}\mj_{B})))\mw_{A{\cdot}B},{\mbox{\hspace{0.3cm}}}
f:A \str C,{\mbox{\hspace{0.3cm}}}g:B \str D
\end{eqnarray*}
and in the standard axiomatization,
\[ \mpm_{A,B}=\mpm_{A,\mi}{\langle}\mj_{A}\mpm_{A,B},\mk_{B}\mpm_{A,B}'
{\rangle},{\mbox{\hspace{0.3cm}}}
\mpm_{A,B}'=\mpm_{\mi,B}{\langle}\mk_{A}\mpm_{A,B},\mj_{B}\mpm_{A,B}'
{\rangle} \]
\[ {\langle}f,g{\rangle}={\langle}f\mpm_{C,C},g\mpm_{C,C}'{\rangle}
{\langle}\mj_{C},\mj_{C}{\rangle},{\mbox{\hspace{0.3cm}}}
f:C \str A,{\mbox{\hspace{0.3cm}}}g:C \str B.\]

\noindent Some derivations are straightforward and some, 
like those concerning
{\bf b}-arrows, require some effort to be proven. We leave them as
an exercise and suggest to the reader who is familiar with the
coherence results in symmetric monoidal categories to use them in these
proofs.

\section{G-natural transformations and transformational graphs}
Let $\cal A$ be an arbitrary category and
$F:{\cal A}^{m}{\rightarrow}
{\cal A}$, $G:{\cal A}^{n}{\rightarrow}{\cal A}$, $m,n \geq 0$, two functors
(${\cal A}^{0}$ is trivial category).
Let $\Gamma$ be a function from $\{1, \ldots , n\}$  to
$\{1, \ldots , m\}$ called {\em graph} (if $n=0$, then
$\{1, \ldots , n\}$ is $\emptyset$). We say that an indexed family of
morphisms from $\cal A$
\[ \alpha{\!}={\!}\{ \alpha(A_{1},\ldots,A_{m}){\!}:{\!}
F(A_{1},\ldots,A_{m})
 \str G(A_{\Gamma(1)},\ldots,A_{\Gamma(n)})\:
{\mid}\:A_{1},{\ldots},A_{m}{\in}{\cal A} \}
\]
is a {\em g-natural transformation} from $F$ to $G$ with the graph $\Gamma$,
denoted by $\alpha:F \gna G$,
if for every $i$, $1{\leq}i{\leq}m$, arbitrary $A_{1},{\ldots},
A_{i},A_{i}',A_{i+1},{\ldots},A_{m}$ and $f:A_{i} \str A_{i}'$
from ${\cal A}$, the following diagram commutes

\begin{center}
\begin{picture}(370,90)
\put(110,10){\vector(1,0){130}}
\put(110,70){\vector(1,0){130}}
\put(50,60){\vector(0,-1){40}}
\put(300,60){\vector(0,-1){40}}
\put(175,14){\makebox(0,0)[b]{$\alpha(A_{1},{\ldots},A_{i}',{\ldots},
A_{m})$}}
\put(175,74){\makebox(0,0)[b]{$\alpha(A_{1},{\ldots},A_{i},{\ldots},
A_{m})$}}
\put(51,40){\makebox(0,0)[l]{$F(\mj_{A_{1}},{\ldots},f,{\ldots},
\mj_{A_{m}})$}}
\put(299,40){\makebox(0,0)[r]{$F(h_{\Gamma(1)},{\ldots},h_{\Gamma(n)})$}}
\put(107,10){\makebox(0,0)[r]{$F(A_{1},{\ldots}, A_{i}', {\ldots},
A_{m})$}}
\put(243,10){\makebox(0,0)[l]{$G(A_{\Gamma(1)}^{\ast},{\ldots},
A_{\Gamma(n)}^{\ast})$}}
\put(107,70){\makebox(0,0)[r]{$F(A_{1},{\ldots},A_{i}, {\ldots},
A_{m})$}}
\put(243,70){\makebox(0,0)[l]{$G(A_{\Gamma(1)},{\ldots},
A_{\Gamma(n)})$}}
\end{picture}
\end{center}

\noindent where for $j{\neq}i$,
$h_{j}{\equiv}\mj_{A_{j}}$, $A_{j}^{\ast}{\equiv}A_{j}$ and for $i$,
$h_{i}{\equiv}f$, $A_{i}^{\ast}{\equiv}A_{i}'$.
\\[0.3cm]
This definition follows the one given in (\cite{mvfc}, page 94).
\\[0.3cm]
{\bf Example}\hspace{2em}
Let $\cal A$ be a relevant category. Denote by \mw\
the indexed set
\[ \{\mw_{A}:A \str A \cdot A \mid A \in {\cal A}\}. \]
Then by the equation (\mw), \mw\ is a 
g-natural transformation from the identity
functor $1_{\cal A}: {\cal A} \str {\cal A}$ to the multiplication
functor $\cdot : {\cal A}^{2} \str {\cal A}$ with the graph
$\Gamma: \{ 1,2 \} \str \{ 1 \}$, $\Gamma(1) = \Gamma(2) = 1$.
\\[0.3cm]
If $\alpha:F{\gna}G$ for $F$, $G$ and $\Gamma$ as above,
then it is easy to see that $\alpha$ is a classical natural transformation
between $F$ and $G':{\cal A}^{m}{\rightarrow}
{\cal A}$ where
\[ G'(A_{1},{\ldots},A_{m})=_{df}G(A_{\Gamma(1)},{\ldots},A_{\Gamma(n)}), \]
\[ G'(f_{1},{\ldots},f_{m})=_{df}G(f_{\Gamma(1)},{\ldots},f_{\Gamma(n)}). \]

\noindent As in the classical case, g-natural transformations can be
composed in the following way. Let $F:{\cal A}^{m}{\rightarrow}{\cal A}$,
$G:{\cal A}^{n}{\rightarrow}{\cal A}$ and
$H:{\cal A}^{l}{\rightarrow}{\cal A}$ be functors.
Let ${\alpha}:F{\fna}G$ and ${\beta}:G{\pna}H$ for some graphs
$\Phi:\{1, \ldots ,n\} \str \{1, \ldots ,m\}$ and
$\Psi:\{1, \ldots ,l\} \str \{1, \ldots ,n\}$. We define
its composition as
\[ {\beta}{\alpha}=_{df}\{ {\beta}{\alpha}(A_{1},{\ldots},
A_{m}){\equiv}{\beta}(A_{\Phi(1)},{\ldots},A_{\Phi(n)})
{\alpha}(A_{1},{\ldots},A_{m}){\mid}
(A_{1},{\ldots},A_{m}){\in}{\cal A}^{m} \}.  \]
Then it is easy to prove (as in the case of classical natural
transformations) that $\beta\alpha$
is a g-natural transformation from $F$ to $H$ with the graph $\Phi \Psi$
(the usual composition of functions  $\Psi$ and $\Phi$).

\noindent Generalization of g-natural transformations to the case of several
categories (here we have only 
$\cal A$) is not essential, and serves just to  complicate
the notation.

\section{Canonical transformations in substructural categories}
Throughout this section, $\cal A$ denotes an arbitrary
substructural category.

\noindent Let $\cal F$ be a set of terms  obtained from symbols $\Box$, \ri\
and binary operation $\cdot$. Its elements we call {\em shapes}.

\noindent In a natural way, we define correspondence between shapes and
functors of type
${\cal A}^{n}{\rightarrow}{\cal A}$ for some $n{\geq}0$.
\begin{enumerate}
\item Functor
$\mj_{\cal A}:{\cal A}{\rightarrow}{\cal A}$ corresponds to the term
$\Box$.
\item Functor
$I:{\cal A}^{0}{\rightarrow}{\cal A}$, which maps the unique object from
${\cal A}^0$
to the object {\ri} from $\cal A$, corresponds to the term \ri.
\item If
$F:{\cal A}^{m}{\rightarrow}{\cal A}$ corresponds to the term
\mfv\ and $G:{\cal A}^{n}{\rightarrow}{\cal A}$ corresponds
to the term \mgv,
then the functor $H:{\cal A}^{m+n}{\rightarrow}{\cal A}$
such that for every $m+n$-tuple $(A_{1},{\ldots},A_{m+n})$
of ${\cal A}$ objects, $H(A_{1},{\ldots},A_{m+n})=_{df}
F(A_{1},{\ldots},A_{m}){\cdot}G(A_{m+1},{\ldots},A_{m+n})$,
and for every $m+n$-tuple $(f_{1},{\ldots},f_{m+n})$ of ${\cal A}$
arrows, $H(f_{1},{\ldots},f_{m+n})=_{df}
F(f_{1},{\ldots},f_{m}){\cdot}G(f_{m+1},{\ldots},f_{m+n})$, corresponds
to the term $\mfv \cdot \mgv$.
\end{enumerate}

\noindent However, depending on the category $\cal A$, two different
shapes may define the same functor. From now on, if we say that a
functor $F:{\cal A}^{n}
{\rightarrow}{\cal A}$ is from $\cal F$, that means it corresponds
to some shape from $\cal F$.

\noindent Functors from $\cal F$ will serve as domains and codomains
of canonical transformations we are going to introduce below.

Let $F:{\cal A}^{m}{\rightarrow}{\cal A}$,
$G:{\cal A}^{n}{\rightarrow}{\cal A}$ and
$H:{\cal A}^{l}{\rightarrow}{\cal A}$ be functors from $\cal F$.
\\[0.2cm]
{\em a)}{\hspace{0.3cm}} Denote by $\mj_{F}$ the indexed family 
$\{ \mj_{F(A_{1},{\ldots},A_{m})}{\mid}(A_{1},{\ldots},A_{m})
{\in}{\cal A}^{m} \}$ and let $\Gamma$ be  the identity function
from $\{1, \ldots , m\}$ to $\{1, \ldots , m\}$. It is easy to
see that $\mj_{F}$ is
a g-natural transformation from $F$ to $F$ with the graph $\Gamma$.
\\[0.3cm]
If ${\cal A}$ is monoidal
\\
{\em b)}{\hspace{0.3cm}}Denote by $\ms_{F}$ the indexed family
$\{ \ms_{F(A_{1},{\ldots},A_{m})}{\mid}(A_{1},{\ldots},A_{m})
{\in}{\cal A}^{m} \}$ and let $\Gamma$ be as above. Then
$\ms_{F}:\ri{\cdot}F{\gna}F$.

\noindent In a similar way we define
$\ms_{F}^{i}$, $\md_{F}$ and
$\md_{F}^{i}$.
\\[0.3cm]
{\em c)}{\hspace{0.3cm}}Denote by $\br_{F,G,H}$ the indexed family
\\
$\{ \br_{F(A_{1},{\ldots},A_{m}),G(A_{m+1},{\ldots},A_{m+n}),
H(A_{m+n+1},{\ldots},A_{m+n+l})
}{\mid}(A_{1},{\ldots},A_{m+n+l})
{\in}{\cal A}^{m+n+l} \}$
and let $\Gamma$ be the identity function from $\{1, \ldots , m+n+l\}$
to  $\{1, \ldots , m+n+l\}$.
Then,
$\br_{F,G,H}:F{\cdot}(G{\cdot}H){\gna}(F{\cdot}G){\cdot}H$.

\noindent In a similar way we define $\bl_{F,G,H}$. 
\\[0.3cm]
If $\cal A$ is symmetric monoidal
\\
{\em d)}{\hspace{0.3cm}}Denote by $\mc_{F,G}$ the indexed set
$\{ \mc_{F(A_{1},{\ldots},A_{m}),G(A_{m+1},{\ldots},A_{m+n})}
{\mid}(A_{1},{\ldots},A_{m+n})
{\in}{\cal A}^{m+n} \}$ and let $\Gamma$ be the function from
$\{1, \ldots , m+n\}$ to $\{1, \ldots , m+n\}$ that satisfies
$\Gamma(m+i)=i$
for $1{\leq}i{\leq}n$ and $\Gamma(j)=n+j$ for $1{\leq}j{\leq}m$.
Then
$\mc_{F,G}:F{\cdot}G{\gna}G{\cdot}F$.
\\[0.3cm]
{\em e)}{\hspace{0.3cm}}If $\cal A$ is relevant category,
we denote by $\mw_{F}$ the indexed family
$\{ \mw_{F(A_{1},{\ldots},A_{m})}
\;{\mid}\;(A_{1},{\ldots},A_{m})
{\in}{\cal A}^{m} \}$. Let $\Gamma$ be the function
from $\{1, \ldots , 2m\}$ to $\{1, \ldots , m\}$ defined as
$\Gamma(i)=\Gamma(m+i)=i$
for $1{\leq}i{\leq}m$. Then
$\mw_{F}:F{\gna}F{\cdot}F$.
\\[0.3cm]
{\em f)}{\hspace{0.3cm}}If $\cal A$ is affine,
we denote by $\mk_{F}$ the indexed family
$\{ \mk_{F(A_{1},{\ldots},A_{m})}
{\mid}(A_{1},{\ldots},A_{m})
{\in}{\cal A}^{m} \}$. Let $\Gamma$ be the empty function
from $\emptyset$ to $\{1, \ldots , m\}$.
Then 
$\mk_{F}:F{\gna}I$.
\\[0.3cm]
The g-natural transformations from above, which exist in the category
$\cal A$ constitute the class of $\cal A$ {\em primitive canonical 
transformations}. If we declare $\cal A$ is a 
monoidal category, its primitive
canonical transformations are those from {\em a)} to {\em c)}, though
$\cal A$ may have the structure of a cartesian category. 
\\[0.1cm]
Let $F_{1}:{\cal A}^{m}{\rightarrow}{\cal A}$,
$F_{2}:{\cal A}^{n}{\rightarrow}{\cal A}$,
$G_{1}:{\cal A}^{k}{\rightarrow}{\cal A}$,
$G_{2}:{\cal A}^{l}{\rightarrow}{\cal A}$ be functors from $\cal F$,
and $\alpha:F_{1}{\fna}G_{1}$ and $\beta:F_{2}{\pna}G_{2}$. Denote by
$\alpha{\cdot}\beta$ the family
\[\{ \alpha(A_{1}, \ldots , A_{m}) \cdot \beta(A_{m+1}, \ldots , A_{m+n})
{\mid} A_{1}, \ldots , A_{m+n} \in {\cal A} \}, \]
and let
$\Gamma$ as a function from $\{1, \ldots , k+l\}$ to
$\{1, \ldots , m+n\}$  satisfy the following:
\\
$\Gamma(i)=\Phi(i)$ for
$1{\leq}i{\leq}k$ and $\Gamma(k+j)=m+\Psi(j)$ for $(1{\leq}j{\leq}l)$.
\noindent Then it is easy to see that
$\alpha{\cdot}\beta:F_{1}{\cdot}F_{2}{\gna}G_{1}{\cdot}G_{2}$.
\\[0.3cm]
Now, we can define {\em canonical transformations} in $\cal A$ as follows
\begin{enumerate}
\item Primitive canonical transformations from $\cal A$ are canonical
transformations.
\item If $\alpha:F_{1}{\fna}G_{1}$ and $\beta:F_{2}{\pna}G_{2}$
are canonical transformations in $\cal A$, then \linebreak
$\alpha{\cdot}\beta:F_{1}{\cdot}F_{2}{\fpna}G_{1}{\cdot}G_{2}$
is canonical.
\item If $\alpha:F{\fna}G$ and $\beta:G{\pna}H$ are canonical
transformations in $\cal A$, then
\linebreak $\beta\alpha:F{\pfna}H$ is canonical.
\end{enumerate}
It is easy to verify that symmetric monoidal canonical transformations
have {\em bijections} as graphs, relevant canonical transformations have
{\em onto} functions as graphs and affine canonical transformations
have {\em one-one} functions as graphs.
\noindent Now we can reformulate MacLane's results from \cite{natass}  in the
following manner:
\\[0.2cm]
{\em If} $\cal A$ {\em is monoidal or symmetric monoidal category
and} $\alpha, \beta: F \gna G$ {\em are two canonical transformations
(with the same graph} $\Gamma${\em ), then} $\alpha$ {\em and} $\beta$
{\em are the same
indexed sets (i.e., the same functions from the sequences of objects
to the morphisms of} $\cal A${\em ).}
\\[0.2cm]
This property of a category we call {\em coherence}.
It completely follows the notion of coherence given in \cite{clcat}.
We can extend this definition to an arbitrary substructural category
$\cal A$. Namely, we say that an arbitrary substructural category
is {\em coherent} if for every pair of canonical transformations
$\alpha$ and $\beta$ of the same type and with the same graph
we have that $\alpha = \beta$ as indexed sets.

\section{Categories {\bf Mon}, {\bf SyMon}, {\bf Rel},
{\bf Aff} and {\bf Cart}}

\noindent Let $\cal M$ be the category whose objects are
monoidal categories
and whose arrows are the monoidal structure preserving functors 
in the language given above. The equational axiomatization of
monoidal categories enables us to distinguish a category from $\cal M$
freely generated by a set of objects. Let $P$ be an infinite
linearly ordered set of
objects, whose elements we call {\em letters}.
We denote by {\bf Mon} the free monoidal category generated by $P$
whose construction is given below. In the same way we will introduce
the categories {\bf SyMon}, {\bf Rel}, {\bf Aff} and {\bf Cart}, namely
the free symmetric monoidal, relevant, affine and cartesian category
generated by the same set $P$ of objects.
\\[0.1cm]
The constructions of these categories are algebraic and the set of objects
is always the set $\cal O$ of terms freely
generated by $P{\cup}\{{\ri}\}$ using the binary operation $\cdot$.
\\[0.1cm]
{\em Primitive morphism-terms} are in the case of {\bf Mon}
\[ \mj_{A}:A \str A \]
\[ \ms_{A}:\ri{\cdot}A \str A {\mbox{\hspace{2cm}}}
\ms_{A}^{i}:A \str \ri{\cdot}A \]
\[ \md_{A}:A{\cdot}\ri \str A {\mbox{\hspace{2cm}}}
\md_{A}^{i}:A \str A{\cdot}\ri \]
\[ \br_{A,B,C}:A{\cdot}(B{\cdot}C) \str (A{\cdot}B){\cdot}C
{\mbox{\hspace{2cm}}} 
\bl_{A,B,C}:(A{\cdot}B){\cdot}C \str A{\cdot}(B{\cdot}C) \]
for all $A,B,C{\in}{\cal O}$.
\\[0.1cm]
{\bf SyMon} {\em primitive morphism-terms} are those of {\bf Mon}
together with
\[ \mc_{A,B}: A \cdot B \str B \cdot A \]
for every $A,B \in {\cal O}$.
\\[0.1cm]
{\bf Rel} {\em primitive morphism-terms} are those of {\bf SyMon}
together with
\[ \mw_{A}: A  \str A \cdot A \]
for every $A \in {\cal O}$.
\\[0.1cm]
{\bf Aff} {\em primitive morphism-terms} are those of {\bf SyMon}
together with
\[ \mk_{A}: A \str \ri \]
for every $A \in {\cal O}$.
\\[0.1cm]
{\bf Rel} and {\bf Aff} primitive morphism-terms make the class
of {\bf Cart} primitive morphism-terms.
\\[0.1cm]
{\em Morphism-terms} are built from the primitive morphism-terms
with the help of the binary operations of composition and
multiplication.
\\[0.1cm]
{\em Morphisms} of the category {\bf Mon} are equivalence classes
of {\bf Mon} morphism-terms modulo monoidal equations. Analogously,
we define morphisms of other
free categories mentioned above.
\\[0.1cm]
Let $\cal C$ be one among
{\bf Mon}, {\bf SyMon}, {\bf Rel},
{\bf Aff} and {\bf Cart}.
We define a correspondence between the morphism-terms 
and canonical transformations of $\cal C$ in the following way
\begin{enumerate}
\item If $f:A \str B$ is primitive morphism-term, suppose it is of the form
$\mj_{F(p_{1},{\ldots},p_{m})}$  for some $F:{\cal C}^{m} \str {\cal C}
{\in}{\cal F}$
and some, not necessarily distinct, letters $p_{1},{\ldots},p_{m}$.
Then  the canonical transformation
$\mj_{F}:F{\gna}F$ where $\Gamma$ is the identity function on
$\{1, \ldots ,m \}$ corresponds to $f$.
We procede similarly in the remaining cases.
\item If $f$ is of the form $f_{1}{\cdot}f_{2}$ or $f_{2}f_{1}$
and if the canonical transformations $\alpha_{1}$ and $\alpha_{2}$
correspond to the morphism-terms $f_{1}$ and $f_{2}$, then 
the canonical transformation $\alpha_{1}{\cdot}\alpha_{2}$ respectively
$\alpha_{2}\alpha_{1}$ corresponds to the morphism-term $f$.
\end{enumerate}
The graph  of a transformation that
corresponds to the morphism-term $f:A \str B$ we call simply a {\em graph}
of $f$. It connects the occurrence of
a letter  in $A$ with a set (maybe empty) of occurrences of the same letter
in $B$.
\\[0.1cm]
Let $\alpha$ be a term of a $\cal C$ canonical transformation of
type $F{\gna}G$,
for $F:{\cal C}^{m} \str {\cal C}$, $G:{\cal C}^{n} \str {\cal C}$
$\in {\cal F}$.
Let $p_{1},p_{2},{\ldots},p_{m}$ be distinct letters from $P$. We call
the morphism-term
\[ \alpha(p_{1},{\ldots},p_{m}):
F(p_{1}, \ldots ,p_{m})
 \str 
G(p_{\Gamma(1)}, \ldots ,p_{\Gamma(n)}) \]
the {\em representative} of the transformation $\alpha$.
\\[0.3cm]
{\bf Lemma 1}\hspace{2em}
{\em
Let $\cal A$ be an arbitrary substructural category and $\cal C$
one of the free categories mentioned above which is
of the same type as $\cal A$. Let $F$ and $G$  be from $\cal F$  and let
$\alpha:F{\fna}G$ and $\beta:F{\pna}G$  be in $\cal A$. Denote by
{\boldmath{$\alpha$}} and {\boldmath{$\beta$}} canonical transformations
in $\cal C$ defined by the same terms as $\alpha$ and $\beta$
respectively.
Let  $f:A \str B$ be the representative of
{\boldmath$\alpha$} and $g:A \str B$  the representative of
{\boldmath$\beta$}. If $f=g$ in $\cal C$, then
$\alpha = \beta$ in $\cal A$.
}
\\[0.3cm]
{\bf proof}\hspace{2em}
Suppose that
$f \equiv \mal (p_{1}, \ldots ,p_{m}):
F(p_{1}, \ldots ,p_{m}) \str G(p_{\Phi(1)}, \ldots ,p_{\Phi(n)})$
for some distinct letters $p_{1}, \ldots , p_{m} \in P$.
Since $g$ is the representative equal to $f$, they share domains 
and codomains; hence
$g \equiv \mbe (p_{1}, \ldots ,p_{m}):
F(p_{1}, \ldots ,p_{m}) \str G(p_{\Psi(1)}, \ldots ,p_{\Psi(n)})$,
which implies that for every $1 \leq i \leq n$, $\Phi(i) = \Psi(i)$, and so
$\Phi$ and $\Psi$ are equal graphs.
Suppose that $f' \equiv \alpha (A_{1}, \ldots ,A_{m})$. By the assumption
that $\cal C$ is free, there is a functor
$U:{\cal C}{\rightarrow}{\cal A}$ that preserves the structure of $\cal C$
and that extends the mapping of the generators given by
$p_{1}{\mapsto}A_{1},\;p_{2}{\mapsto}A_{2},{\ldots},p_{m}{\mapsto}A_{m}$
(other generators are mapped arbitrarily). Then we have
\[ \alpha{\ni}f'={\alpha}(A_{1},{\ldots},A_{k})=U(f)=U(g)=
{\beta}(A_{1},{\ldots},A_{k}){\in}\beta, \]
hence $\alpha{\subset}\beta$. In the same way we prove
that $\beta{\subset}\alpha$.
\qed
Our goal is to prove that every substructural category is coherent,
and the following lemma will serve to reduce the problem to the case
of the free category in the type. For an object $A \in {\cal O}$ we say
that it is {\em diversified} if no letter occurs twice in it.
\\[0.3cm]
{\bf Lemma 2}\hspace{2em}
{\em
If for every pair of $\cal C$ morphism-terms $f,g:A \str B$,
such that $A$ is diversified, holds that $f=g$,
then every category of the same type as $\cal C$ is coherent.
}
\\[0.3cm]
{\bf proof}\hspace{2em}
Let $\cal A$ be an arbitrary substructural category of the same type as
$\cal C$. Suppose that $F$ and $G$ are functors from $\cal F$ and
$\alpha, \beta : F \gna G$ in $\cal A$.
Let $f \equiv \mal(p_{1}, \ldots , p_{m}) : A \str B$
and $g \equiv \mbe(p_{1}, \ldots , p_{m}) : C \str D$
be the representatives of \mal\ and \mbe\ respectively, where
\mal\ and \mbe\ are $\cal C$ canonical transformations 
defined by the same terms as $\alpha$ and $\beta$.
Since they have the same graph, $A$ must be identical to $C$ and
$B$ to $D$. By the definition
of graph, it follows that $A$ is diversified
and, by the assumption, $f$ is equal to $g$.
Hence, by Lemma 1, we have $\alpha = \beta$.
\qed
By the following series of definitions we introduce some auxiliary
notions that will help us in proving our coherence results.

Denote by $\cal PF$ the set of terms generated by the binary operation
$\cdot$ from the elements of $P \cup \{ \Box, \ri \}$
(e.g. $(\Box \cdot p) \cdot ((\ri \cdot \Box) \cdot q)$ is in
$\cal PF$).  As in the case of terms from $\cal F$,
we can define in the same way the correspondence
between terms from $\cal PF$
and  functors (with parameters) of the type ${\cal C}^{n} \str {\cal C}$
for some $n \geq 0$, where $\cal C$ is one of the free categories
mentioned above.
\\[0.3cm]
A {\em product term} of $\cal C$ is a morphism-term
defined recursively as follows
\begin{enumerate}
\item The primitive terms (if they exist in $\cal C$)
\[ \ms_{Q},\ms^{i}_{Q},\md_{Q},\md^{i}_{Q},\br_{Q,S,R},
\bl_{Q,S,R},\mc_{Q,S},\mw_{Q},\mk_{Q}. \]
are product terms, called {\em determining factors}.
\item The terms $\mj_{Q}$ are product terms.
\item If $f$ is a product term, then $\mj_{Q}{\cdot}f$ and
$f{\cdot}\mj_{Q}$ are product terms.
\end{enumerate}
	The determining factor of a product term $f$, if it exists, is
denoted  by $d(f)$ (we call such a term {\em structural product}).
A structural product $f$ is a $\mb$-product iff $d(f)$ is a
$\mb$
term, $\mc$-product iff $d(f)$ is a $\mc$ term, and similarly for 
$\ms$,
$\ms ^{i}$, $\md$, $\md ^{i}$, $\mk$ and $\mw$-products.
\\[0.3cm]
For a \mw-product-term we say that it is {\em atomic}
if the index of its determining factor is a letter.
\\[0.3cm]
We say that an atomic \mw-product is {\em left} if there is not
any \mj\ with the letter $p$ in the index, on the left of
its determining factor
$\mw_{p}$ (e.g.
$(\mj_{q{\cdot}r}{\cdot}\mw_{p}){\cdot}\mj_{p}$
is the left atomic \mw-product, while
$(\mj_{p{\cdot}r}{\cdot}\mw_{p}){\cdot}\mj_{p}$
is not left.)
\\[0.3cm]
For a \mc-product we say that it is {\em atomic} if the index of
its determining factor
is a pair of atoms (an atom is a letter or  \ri).
\\[0.3cm]
We say that an atomic \mc-product is {\em diversified} if its
determining factor is not of the form $\mc_{p,p}$ for some letter $p$.
\\[0.3cm]
For a \mk-product we say that it is {\em atomic} if the index
of its determining factor is a letter.
\\[0.3cm]
We say that a composition of atomic \mw-products  (\mk-products) 
is {\em ordered}, if a $\mw_{p}$-product ($\mk_{p}$-product)
is to the right of $\mw_{q}$-product ($\mk_{q}$-product) in this
composition
iff the letter $p$ precedes the letter $q$ in the ordering of $P$.
\section{Coherence in Relevant categories}
At the beginning of this section, we prove a lemma that states 
coherence for {\bf bw} fragments of relevant categories. 
\\[0.3cm]
{\bf Lemma 3}\hspace{2em}
{\em
Let $F$ be from $\cal F$ and let $f:p{\:} \str F(p, \ldots ,p)$, be
a composition of atomic \mw-products.
Then $f$ is equal to a term of the form $hg$ where $g$
is a composition of left atomic \mw-products and $h$ is a
composition of {\bf b}-products.
}
\dkz
For the sake of clarity we introduce a tree that corresponds to $f$,
denoted by $\tau_{f}$, in the following way.

If $f{\equiv}\mw_{p}$, then $\tau_{f}$ is
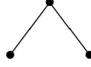
\begin{center}
\begin{picture}(30,20)
\put(0,0){\circle*{3}}
\put(30,0){\circle*{3}}
\put(15,20){\circle*{3}}
\put(0,0){\line(3,4){15}}
\put(15,20){\line(3,-4){15}}
\end{picture}
\end{center}

If $f$ is of the form $G(\mw_{p})\:f_{1}$, where $G$ is from $\cal PF$,
and if in the shape of $G$, $i-1$ letters $p$ precede (from the left)
the symbol $\Box$, then $\tau_{f}$ is obtained
from $\tau_{f_{1}}$ by forking the $i$-th leaf
(from the left) and concatenating simple segments to remaining leaves.
\\
For example, if $f$ is of the form
\[ ((\mj_{p}{\cdot}\mw_{p}){\cdot}\mj_{(p{\cdot}p){\cdot}p})
(\mj_{p{\cdot}p}{\cdot}(\mw_{p}{\cdot}\mj_{p}))
(\mj_{p{\cdot}p}{\cdot}\mw_{p})(\mw_{p}{\cdot}\mj_{p})\mw_{p} \]
then the corresponding tree is
\begin{center}
\begin{picture}(60,50)
\put(60,0){\circle*{3}}
\put(45,0){\circle*{3}}
\put(35,0){\circle*{3}}
\put(25,0){\circle*{3}}
\put(15,0){\circle*{3}}
\put(0,0){\circle*{3}}
\put(60,10){\circle*{3}}
\put(45,10){\circle*{3}}
\put(35,10){\circle*{3}}
\put(20,10){\circle*{3}}
\put(0,10){\circle*{3}}
\put(60,20){\circle*{3}}
\put(40,20){\circle*{3}}
\put(20,20){\circle*{3}}
\put(0,20){\circle*{3}}
\put(50,30){\circle*{3}}
\put(20,30){\circle*{3}}
\put(0,30){\circle*{3}}
\put(50,40){\circle*{3}}
\put(10,40){\circle*{3}}
\put(30,50){\circle*{3}}
\put(0,30){\line(0,-1){30}}
\put(15,0){\line(1,2){5}}
\put(25,0){\line(-1,2){5}}
\put(35,0){\line(0,1){10}}
\put(45,0){\line(0,1){10}}
\put(60,0){\line(0,1){20}}
\put(20,10){\line(0,1){20}}
\put(35,10){\line(1,2){5}}
\put(45,10){\line(-1,2){5}}
\put(40,20){\line(1,1){10}}
\put(60,20){\line(-1,1){10}}
\put(0,30){\line(1,1){10}}
\put(20,30){\line(-1,1){10}}
\put(50,30){\line(0,1){10}}
\put(10,40){\line(2,1){20}}
\put(50,40){\line(-2,1){20}}
\end{picture}
\end{center}
Denote by  $\Lambda$ the set of forking vertices in such a tree.
Let $k_{\lambda}$ be the number of right branches
of these forkings in the path from the vertex $\lambda$ to the root.
The complexity of the tree is measured by the number $n_{f}$ defined as
\[ n_{f}=_{df}\sum_{{\lambda}{\in}{\Lambda}} k_{\lambda} \]
In the example above $n_{f}$ is $3$.
We prove the lemma by induction on $n_{f}$.

If $n_{f}=0$, then $f$ is itself a composition of left atomic
\mw-products.

If $n_{f}>0$ and if there is no subtree of $\tau_{f}$
of the form
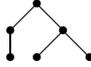
\begin{center}
\begin{picture}(30,20)
\put(10,20){\circle*{3}}
\put(0,10){\circle*{3}}
\put(20,10){\circle*{3}}
\put(0,0){\circle*{3}}
\put(10,0){\circle*{3}}
\put(30,0){\circle*{3}}
\put(0,0){\line(0,1){10}}
\put(10,0){\line(1,1){10}}
\put(30,0){\line(-1,1){20}}
\put(0,10){\line(1,1){10}}
\end{picture}
\end{center}
then by the functoriality of $\cdot$
we obtain a term $f'$ equal to $f$ such that
$n_{f'}=n_{f}$, and there is a subtree of the above form in $\tau_{f'}$. 
In the example,
we obtain the term
\[ ((\mj_{p}{\cdot}\mw_{p}){\cdot}\mj_{(p{\cdot}p){\cdot}p})
(\mj_{p{\cdot}p}{\cdot}(\mw_{p}{\cdot}\mj_{p}))
(\mw_{p}{\cdot}\mj_{p{\cdot}p})(\mj_{p}{\cdot}\mw_{p})\mw_{p}, \]
whose tree is
\begin{center}
\begin{picture}(80,50)
\put(60,0){\circle*{3}}
\put(45,0){\circle*{3}}
\put(35,0){\circle*{3}}
\put(25,0){\circle*{3}}
\put(15,0){\circle*{3}}
\put(0,0){\circle*{3}}
\put(60,10){\circle*{3}}
\put(45,10){\circle*{3}}
\put(35,10){\circle*{3}}
\put(20,10){\circle*{3}}
\put(0,10){\circle*{3}}
\put(60,20){\circle*{3}}
\put(40,20){\circle*{3}}
\put(20,20){\circle*{3}}
\put(0,20){\circle*{3}}
\put(60,30){\circle*{3}}
\put(40,30){\circle*{3}}
\put(10,30){\circle*{3}}
\put(50,40){\circle*{3}}
\put(10,40){\circle*{3}}
\put(30,50){\circle*{3}}
\put(0,0){\line(0,1){20}}
\put(15,0){\line(1,2){5}}
\put(25,0){\line(-1,2){5}}
\put(35,0){\line(0,1){10}}
\put(45,0){\line(0,1){10}}
\put(60,0){\line(0,1){30}}
\put(20,10){\line(0,1){10}}
\put(35,10){\line(1,2){5}}
\put(45,10){\line(-1,2){5}}
\put(40,20){\line(0,1){10}}
\put(0,20){\line(1,1){10}}
\put(20,20){\line(-1,1){10}}
\put(10,30){\line(0,1){10}}
\put(40,30){\line(1,1){10}}
\put(60,30){\line(-1,1){10}}
\put(50,40){\line(-2,1){20}}
\put(10,40){\line(2,1){20}}
\put(80,20){\shortstack{$n=3$}}
\end{picture}
\end{center}

By using the ({\bf bw}) equality and the naturality of {\bf b}-products
we transform this term into the form $h_{1}f_{1}$ where $f_{1}$
is a composition of atomic \mw-products with
$n_{f'}=n_{f}-1$, and $h_{1}$ is a composition of {\bf b}-products.
In our example we obtain the term
\[ \bl_{p{\cdot}(p{\cdot}p),p{\cdot}p,p}
(((\mj_{p}{\cdot}\mw_{p}){\cdot}\mj_{p{\cdot}p}){\cdot}\mj_{p})
((\mj_{p{\cdot}p}{\cdot}\mw_{p}){\cdot}\mj_{p})
((\mw_{p}{\cdot}\mj_{p}){\cdot}\mj_{p})(\mw_{p}{\cdot}\mj_{p})\mw_{p}, \]
and the tree corresponding to its initial part ($f_{1}$) is
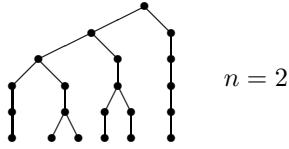
\begin{center}
\begin{picture}(80,50)
\put(60,0){\circle*{3}}
\put(45,0){\circle*{3}}
\put(35,0){\circle*{3}}
\put(25,0){\circle*{3}}
\put(15,0){\circle*{3}}
\put(0,0){\circle*{3}}
\put(60,10){\circle*{3}}
\put(45,10){\circle*{3}}
\put(35,10){\circle*{3}}
\put(20,10){\circle*{3}}
\put(0,10){\circle*{3}}
\put(60,20){\circle*{3}}
\put(40,20){\circle*{3}}
\put(20,20){\circle*{3}}
\put(0,20){\circle*{3}}
\put(60,30){\circle*{3}}
\put(40,30){\circle*{3}}
\put(10,30){\circle*{3}}
\put(60,40){\circle*{3}}
\put(30,40){\circle*{3}}
\put(50,50){\circle*{3}}
\put(0,0){\line(0,1){20}}
\put(15,0){\line(1,2){5}}
\put(25,0){\line(-1,2){5}}
\put(35,0){\line(0,1){10}}
\put(45,0){\line(0,1){10}}
\put(60,0){\line(0,1){40}}
\put(20,10){\line(0,1){10}}
\put(35,10){\line(1,2){5}}
\put(45,10){\line(-1,2){5}}
\put(40,20){\line(0,1){10}}
\put(0,20){\line(1,1){10}}
\put(20,20){\line(-1,1){10}}
\put(10,30){\line(2,1){40}}
\put(40,30){\line(-1,1){10}}
\put(60,40){\line(-1,1){10}}
\put(80,20){\shortstack{$n=2$}}
\end{picture}
\end{center}
By the induction hypothesis the term $f_{1}$ is equal to the term
$h_{2}g$ of the desired form, and therefore $f=h_{1}h_{2}g$ is such.
\qed
In  our example,  the last term is transformed into
\[ \bl_{p{\cdot}(p{\cdot}p),p{\cdot}p,p}
(((\mj_{p}{\cdot}\mw_{p}){\cdot}\mj_{p{\cdot}p}){\cdot}\mj_{p})
((\mw_{p}{\cdot}\mj_{p{\cdot}p}){\cdot}\mj_{p})
((\mj_{p}{\cdot}\mw_{p}){\cdot}\mj_{p})(\mw_{p}{\cdot}\mj_{p})\mw_{p}, \]
and the tree corresponding to its initial part is
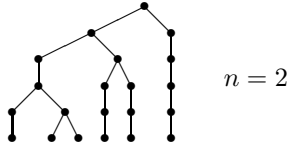
\begin{center}
\begin{picture}(80,50)
\put(60,0){\circle*{3}}
\put(45,0){\circle*{3}}
\put(35,0){\circle*{3}}
\put(25,0){\circle*{3}}
\put(15,0){\circle*{3}}
\put(0,0){\circle*{3}}
\put(60,10){\circle*{3}}
\put(45,10){\circle*{3}}
\put(35,10){\circle*{3}}
\put(20,10){\circle*{3}}
\put(0,10){\circle*{3}}
\put(60,20){\circle*{3}}
\put(45,20){\circle*{3}}
\put(35,20){\circle*{3}}
\put(10,20){\circle*{3}}
\put(60,30){\circle*{3}}
\put(40,30){\circle*{3}}
\put(10,30){\circle*{3}}
\put(60,40){\circle*{3}}
\put(30,40){\circle*{3}}
\put(50,50){\circle*{3}}
\put(0,0){\line(0,1){10}}
\put(15,0){\line(1,2){5}}
\put(25,0){\line(-1,2){5}}
\put(35,0){\line(0,1){20}}
\put(45,0){\line(0,1){20}}
\put(60,0){\line(0,1){40}}
\put(0,10){\line(1,1){10}}
\put(20,10){\line(-1,1){10}}
\put(45,20){\line(-1,2){5}}
\put(35,20){\line(1,2){5}}
\put(10,20){\line(0,1){10}}
\put(10,30){\line(2,1){40}}
\put(40,30){\line(-1,1){10}}
\put(60,40){\line(-1,1){10}}
\put(80,20){\shortstack{$n=2$}}
\end{picture}
\end{center}
(we do this to obtain a subtree of the form
\begin{picture}(35,20)(-5,7)
\put(10,20){\circle*{3}}
\put(0,10){\circle*{3}}
\put(20,10){\circle*{3}}
\put(0,0){\circle*{3}}
\put(10,0){\circle*{3}}
\put(30,0){\circle*{3}}
\put(0,0){\line(0,1){10}}
\put(10,0){\line(1,1){10}}
\put(30,0){\line(-1,1){20}}
\put(0,10){\line(1,1){10}}
\end{picture} )
\\[0.5cm]
and then by ({\bf bw}) and ({\bf b}) this term is transformed into
\[
\bl_{p{\cdot}(p{\cdot}p),p{\cdot}p,p}
(\bl_{p{\cdot}(p{\cdot}p),p,p}{\cdot}\mj_{p}) 
((((\mj_{p}{\cdot}\mw_{p}){\cdot}\mj_{p}){\cdot}\mj_{p}){\cdot}\mj_{p})
((\mw_{p}{\cdot}\mj_{p}){\cdot}\mj_{p}){\cdot}\mj_{p})
((\mw_{p}{\cdot}\mj_{p}){\cdot}\mj_{p})(\mw_{p}{\cdot}\mj_{p})\mw_{p},
\]
to whose initial part corresponds the tree
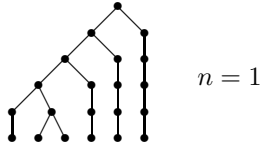
\begin{center}
\begin{picture}(80,50)
\put(50,0){\circle*{3}}
\put(40,0){\circle*{3}}
\put(30,0){\circle*{3}}
\put(20,0){\circle*{3}}
\put(10,0){\circle*{3}}
\put(0,0){\circle*{3}}
\put(30,10){\circle*{3}}
\put(50,10){\circle*{3}}
\put(40,10){\circle*{3}}
\put(15,10){\circle*{3}}
\put(0,10){\circle*{3}}
\put(50,20){\circle*{3}}
\put(40,20){\circle*{3}}
\put(30,20){\circle*{3}}
\put(10,20){\circle*{3}}
\put(50,30){\circle*{3}}
\put(40,30){\circle*{3}}
\put(20,30){\circle*{3}}
\put(50,40){\circle*{3}}
\put(30,40){\circle*{3}}
\put(40,50){\circle*{3}}
\put(0,0){\line(0,1){10}}
\put(10,0){\line(1,2){5}}
\put(20,0){\line(-1,2){10}}
\put(30,0){\line(0,1){20}}
\put(40,0){\line(0,1){30}}
\put(50,0){\line(0,1){40}}
\put(0,10){\line(1,1){40}}
\put(30,20){\line(-1,1){10}}
\put(40,30){\line(-1,1){10}}
\put(50,40){\line(-1,1){10}}
\put(70,20){\shortstack{$n=1$}}
\end{picture}
\end{center}
Then, again by ({\bf bw}) and ({\bf b}) the last term is transformed into
the term
\begin{eqnarray*}
\lefteqn{\bl_{p{\cdot}(p{\cdot}p),p{\cdot}p,p}
(\bl_{p{\cdot}(p{\cdot}p),p,p}{\cdot}\mj_{p})
(((\bl_{p,p,p}{\cdot}\mj_{p}){\cdot}\mj_{p}){\cdot}\mj_{p})} \\
& & ((((\mw_{p}{\cdot}\mj_{p}){\cdot}\mj_{p}){\cdot}\mj_{p}){\cdot}\mj_{p})
(((\mw_{p}{\cdot}\mj_{p}){\cdot}\mj_{p}){\cdot}\mj_{p})
((\mw_{p}{\cdot}\mj_{p}){\cdot}\mj_{p})(\mw_{p}{\cdot}\mj_{p})\mw_{p},
\end{eqnarray*}
of the desired form, whose tree is of the form
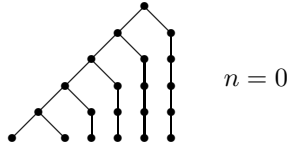
\begin{center}
\begin{picture}(80,50)
\put(50,0){\circle*{3}}
\put(40,0){\circle*{3}}
\put(30,0){\circle*{3}}
\put(20,0){\circle*{3}}
\put(60,0){\circle*{3}}
\put(0,0){\circle*{3}}
\put(30,10){\circle*{3}}
\put(50,10){\circle*{3}}
\put(40,10){\circle*{3}}
\put(60,10){\circle*{3}}
\put(10,10){\circle*{3}}
\put(50,20){\circle*{3}}
\put(40,20){\circle*{3}}
\put(60,20){\circle*{3}}
\put(20,20){\circle*{3}}
\put(50,30){\circle*{3}}
\put(60,30){\circle*{3}}
\put(30,30){\circle*{3}}
\put(60,40){\circle*{3}}
\put(40,40){\circle*{3}}
\put(50,50){\circle*{3}}
\put(0,0){\line(1,1){50}}
\put(20,0){\line(-1,1){10}}
\put(30,0){\line(0,1){10}}
\put(40,0){\line(0,1){20}}
\put(50,0){\line(0,1){30}}
\put(60,0){\line(0,1){40}}
\put(30,10){\line(-1,1){10}}
\put(40,20){\line(-1,1){10}}
\put(50,30){\line(-1,1){10}}
\put(60,40){\line(-1,1){10}}
\put(80,20){\shortstack{$n=0$}}
\end{picture}
\end{center}
{\bf Corollary}{\hspace{1em}}
{\em
Let $F:{\mbox{\bf Rel}}^{\:k}{\rightarrow}{\mbox{\bf Rel}}$ be from
$\cal F$ and let $f:p\: \str F(p,{\ldots},p)$, be
a composition of atomic \mw-products. Then for every $i$, 
$1{\leq}i{\leq}k-1$ there is a morphism term of the form
\[ v((\mj_{\underbrace{p{\cdot}p{\cdot}{\ldots}{\cdot}p}_{i-1}}
{\cdot}\mw_{p}){\cdot}
\mj_{\underbrace{p{\cdot}p{\cdot}{\ldots}{\cdot}p}_{k-i-1}})u \]
equal to $f$ (all products of $p$'s are associated to the left, i.e., 
$p \cdot p \cdot p$ means $(p \cdot p) \cdot p$), where
$u:p\: \str (p{\cdot}p{\cdot}{\ldots}{\cdot}p)
{\cdot} (p{\cdot}p{\cdot}{\ldots}{\cdot}p)$ is a composition of
atomic \mw-products, and $v$ is a composition of {\bf b}-products.}
\\[0.3cm]
\dkz
By Lemma 3, there is a term of the form
$h_{1}g$ equal to $f$, where $g$ is a composition of left atomic
\mw-products, and $h_{1}$ is a composition of {\bf b}-products.
Let
\[
u:p\: \str {\underbrace{(p{\cdot}p{\cdot}{\ldots}{\cdot}p)}_{i}}
{\cdot}{\underbrace{(p{\cdot}p{\cdot}{\ldots}{\cdot}p)}_{k-i-1}}
\]
be a composition of atomic \mw-products (such always exists). Again, by
Lemma 3 there is a term of the form $h_{2}g$ equal to the term
\[ ((\mj_{\underbrace{p{\cdot}p{\cdot}{\ldots}{\cdot}p}_{i-1}}
{\cdot}\mw_{p}){\cdot}
\mj_{\underbrace{p{\cdot}p{\cdot}{\ldots}{\cdot}p}_{k-i-1}})u, \]
where $g$ is as before, and $h_{2}$ is a composition
of {\bf b}-products. Then,
\[ f=h_{1}g=h_{1}h_{2}^{-1}h_{2}g=h_{1}h_{2}^{-1}
((\mj_{(p{\cdot}p){\cdot}{\ldots}{\cdot}p}
{\cdot}\mw_{p}){\cdot}
\mj_{(p{\cdot}p){\cdot}{\ldots}{\cdot}p})u, \]
where $h_{2}^{-1}$ denotes a composition of {\bf b}-products
inverse to  $h_{2}$ (by ({\bf bb}) equalities).
\qed

This corollary will be of use in the proof of the main result of this
section, which states that
\\[0.3cm]
{\bf Theorem 1}\hspace{2em}
{\em  Every relevant category is coherent.}
\\[0.3cm]
In the proof of the theorem, the following lemma, that gives
the normal form of a morphism from {\bf Rel}, is crucial.
\\[0.3cm]
{\bf Lemma 4}\hspace{2em}
{\em
Let $h:A \str B$ be a morphism-term from {\bf Rel} with $A$ diversified.
Then $h$ is equal to the morphism-term
of the form $h''h'$ where $h'$ is an ordered composition of atomic left
\mw-products, and $h''$ is a composition of products with no \mw-products
and with all \mc-products atomic diversified.
}
\dkz
The transformation of the term $h$ is made in several steps.
For the sake of clarity, we illustrate every step starting with the
term
\[ (\mj_{q}{\cdot}(\mw_{p{\cdot}p}\mw_{p}))\mc_{p,q}. \]
$1^{\circ}${\hspace{0.5cm}}
The term $h$ is equal to  a  composition of product terms.
This follows from the functoriality of multiplication.
In our example,
\[ h=(\mj_{q}{\cdot}\mw_{p{\cdot}p})(\mj_{q}{\cdot}\mw_{p})\mc_{p,q}. \]
$2^{\circ}${\hspace{0.5cm}}
By ({\bf bcw}8) and (\ms\md\mw),  $h$  is equal to a term
with all \mw-products  atomic. In our example
\[ h=t(\mj_{q}{\cdot}(\mw_{p}{\cdot}\mj_{p{\cdot}p}))
(\mj_{q}{\cdot}(\mj_{p}{\cdot}\mw_{p}))(\mj_{q}{\cdot}\mw_{p})
\mc_{p,q}, \]
where
\[ t{\equiv}(\mj_{q}{\cdot}\bl_{p{\cdot}p,p,p})(\mj_{q}{\cdot}
(\br_{p,p,p}{\cdot}\mj_{p}))(\mj_{q}{\cdot}
((\mj_{p}{\cdot}\mc_{p,p}){\cdot}\mj_{p}))
(\mj_{q}{\cdot}(\bl_{p,p,p}{\cdot}\mj_{p}))(\mj_{q}{\cdot}\br{p{\cdot}p
,p,p}) \]
$3^{\circ}${\hspace{0.5cm}}
By the multiplication functoriality and the naturality of \ms, \md, {\bf b},
\mc-products,  atomic \mw-products permute towards the right end,
in order to obtain a term whose initial part (from the right) consists
of atomic \mw-products and
whose tail is a \mw-free composition of products. In our example
this term is
\[ t\mc_{(p{\cdot}p){\cdot}(p{\cdot}p),q}
((\mw_{p}{\cdot}\mj_{p{\cdot}p}){\cdot}\mj_{q})((\mj_{p}{\cdot}\mw_{p})
{\cdot}\mj_{q})(\mw_{p}{\cdot}\mj_{q}) \]
$4^{\circ}${\hspace{0.5cm}}
Using ({\bf b}\mc6), this term is equal
to a term whose \mc-products are atomized.
In the example, this atomization is not essential because its
application to the unique
nonatomic
\mc-product $\mc_{(p{\cdot}p){\cdot}(p{\cdot}p),q}$ 
produces only diversified atomic \mc-products (see the following step),
and therefore we write it in abbreviated form, which is enough
for the further analysis:
\[ h=t_{2}(\mj_{q}{\cdot}((\mj_{p}{\cdot}\mc_{p,p}){\cdot}\mj_{p}))
t_{1}((\mw_{p}{\cdot}\mj_{p{\cdot}p}){\cdot}\mj_{q})
((\mj_{p}{\cdot}\mw_{p}){\cdot}\mj_{q})(\mw_{p}{\cdot}\mj_{q}), \]
\[ {\mbox{where}}\; t_{2}{\equiv}
(\mj_{q}{\cdot}\bl_{p{\cdot}p,p,p})(\mj_{q}{\cdot}
(\br_{p,p,p}{\cdot}\mj_{p})), \]
\[ {\mbox{and}}\; t_{1}{\equiv}
(\mj_{q}{\cdot}(\bl_{p,p,p}{\cdot}\mj_{p}))(\mj_{q}{\cdot}\br{p{\cdot}p
,p,p}){(\mc_{(p{\cdot}p){\cdot}(p{\cdot}p),q})}^{\ast} \]
($\ast$ means a developed form with atomic \mc-products.)
\\[0.3cm]
$5^{\circ}${\hspace{0.5cm}}
Suppose that the tail (\mw-free composition of products)
from the last term is in the form $t_{2}F(\mc_{p,p})t_{1}$, where
$F:{\mbox{\bf Rel}}{\rightarrow}{\mbox{\bf Rel}}$
is from $\cal PF$ and all the \mc-products in $t_{1}$ are atomic
diversified.
Let $F(\mc_{p,p})t_{1}$ be of the type $G(p,p) \str F(p{\cdot}p)$ for
some $G:{\mbox{\bf Rel}}^{2}{\rightarrow}{\mbox{\bf Rel}}$
from $\cal PF$, where the left $p$ from $F(p,p)$ is mapped to the right
$p$ from $G(p{\cdot}p)$ by the graph of $F(\mc_{p,p})t_{1}$,
and the right $p$ from $F(p,p)$ is mapped to the left $p$ from
$F(p{\cdot}p)$ by the same graph.
By the assumption that $A$ is diversified,
that all \mc-products in $t_{1}$ are atomic diversified and
the initial part consists of atomic \mw-products only,
we have that two emphasized letters $p$ in $G(p,p)$
occur consecutively (not necessarily in the form $(p{\cdot}p)$).
\\[0.2cm]
Using the functoriality of $\cdot$, we can push all $\mw_{p}$
products to the (left) end of an initial part consisting of \mw-products
only. Now, we apply the corollary of Lemma 3, to show that such an
initial part is equal to a term of the form
$i_{2}H(\mw_{p})i_{1}$, where $i_{1}$ is a composition of atomic
\mw-products, $H:{\mbox{\bf Rel}}{\rightarrow}{\mbox{\bf Rel}}$
is from $\cal PF$,
$i_{2}$ is a composition of {\bf b}-products, the term
$i_{2}H(\mw_{p})$ is of the type $H(p) \str G(p,p)$ and its
graph maps both distinguished $p$'s from
$G(p,p)$ to the distinguished $p$ in $H(p)$.
(Assume that distinguished $(p,p)$ in
$G(p,p)$ are $i$-th and $i+1$-th occurrence of the letter $p$ and just
apply the corollary of Lemma 3.)
\\
In our example
\[ h=t_{2}(\mj_{q}{\cdot}((\mj_{p}{\cdot}\mc_{p,p}){\cdot}\mj_{p}))
t_{1}i_{2}(((\mj_{p}{\cdot}\mw_{p}){\cdot}\mj_{p}){\cdot}\mj_{q})i_{1}, \]
where $i_{1}{\equiv}((\mw_{p}{\cdot}\mj_{p}){\cdot}\mj_{q})
(\mw_{p}{\cdot}\mj_{q})$, $i_{2}{\equiv}(\bl_{p{\cdot}p,p,p}{\cdot}\mj_{q})
((\br_{p,p,p}{\cdot}\mj_{p}){\cdot}\mj_{q})$,
\\
and $F{\equiv}(q{\cdot}(p{\cdot}{\Box})){\cdot}p)$,
$G{\equiv}((p{\cdot}{\Box}){\cdot}({\Box}{\cdot}p)){\cdot}q$,
$H{\equiv}((p{\cdot}{\Box}){\cdot}p){\cdot}q$.
\\[0.3cm]
The {\bf SyMon} terms $F(\mc_{p,p})t_{1}i_{2}$ and $t_{1}i_{2}H(\mc_{p,p})$
are of the same type. Denote by $\alpha$ and $\beta$
canonical transformations in {\bf SyMon} corresponding to these terms.
It is easy to see that
$\alpha$ and $\beta$ have the same graphs (the graph of $H(\mc_{p,p})$
``commutes'' in composition with the graph of $t_{1}i_{2}$ and is transformed
into the graph of $F(\mc_{p,p})$). By  MacLane's coherence 
for symmetric monoidal categories, we have that $\alpha=\beta$.
From the property that every canonical transformation of {\bf SyMon}
contains at most one morphism of a certain type (this is because
the set $\cal O$ of its objects is freely generated by $P \cup \{ \ri \}$),
we conclude that $F(\mc_{p,p})t_{1}i_{2}=t_{1}i_{2}H(\mc_{p,p})$ holds
in {\bf SyMon}. Now, because all the {\bf SyMon}-equalities hold
in {\bf Rel}, these terms are equal in {\bf Rel} too.
\\
Therefore, $h$ is equal to a term of
the form $t_{2}t_{1}i_{2}H(\mc_{p,p})H(\mw_{p})i_{1}$, which is
by (\mc\mw), equal to $t_{2}t_{1}i_{2}H(\mw_{p})i_{1}$.
Repeating this procedure we can eliminate non diversified \mc\
products occurring in $t_{2}$ obtaining a term equal to $h$
whose initial part consists of atomic
\mw-products and whose tail is a \mw-free composition  of products,
whose \mc-products are atomic diversified.
\\
In the example, this procedure includes the following steps
\begin{eqnarray*}
\lefteqn{t_{2}(\mj_{q}{\cdot}((\mj_{p}{\cdot}\mc_{p,p}){\cdot}\mj_{p}))
t_{1}i_{2}(((\mj_{p}{\cdot}\mw_{p}){\cdot}\mj_{p}){\cdot}\mj_{q})i_{1}=}\\
& & t_{2}t_{1}i_{2}
(((\mj_{p}{\cdot}\mc_{p,p}){\cdot}\mj_{p}){\cdot}\mj_{q})
(((\mj_{p}{\cdot}\mw_{p}){\cdot}\mj_{p}){\cdot}\mj_{q})i_{1}=\\
& & t_{2}t_{1}i_{2}(((\mj_{p}{\cdot}\mw_{p})
{\cdot}\mj_{p}){\cdot}\mj_{q})i_{1},
\end{eqnarray*}
where $t_{2}{\equiv}(\mj_{q}{\cdot}\bl_{p{\cdot}p,p,p})
(\mj_{q}{\cdot}(\br_{p,p,p}{\cdot}\mj_{p}))$. 
\\[0.3cm]
$6^{\circ}${\hspace{0.5cm}}
By Lemma 3 and the functoriality of multiplication,
this term  is equal to the one whose
initial part is an ordered composition of
atomic left \mw-products.
\qed
In the example, we have
\begin{eqnarray*}
\lefteqn{(((\mj_{p}{\cdot}\mw_{p}){\cdot}\mj_{p}){\cdot}\mj_{q})i_{1}
{\equiv}}\\
& & (((\mj_{p}{\cdot}\mw_{p}){\cdot}\mj_{p}){\cdot}\mj_{q})
((\mw_{p}{\cdot}\mj_{p}){\cdot}\mj_{q})(\mw_{p}{\cdot}\mj_{q})=\\
& & ((\bl_{p,p,p}{\cdot}\mj_{p}){\cdot}\mj_{q})(((
\mw_{p}{\cdot}\mj_{p}){\cdot}\mj_{p}){\cdot}\mj_{q})
((\mw_{p}{\cdot}\mj_{p}){\cdot}\mj_{q})(\mw_{p}{\cdot}\mj_{q}).
\end{eqnarray*}
\\[0.3cm]
{\bf Lemma 5}\hspace{2em}
{\em
Let $f,g:A \str B$ be two morphism-terms in {\bf Rel} with
$A$ diversified. Then $f=g$ in {\bf Rel}.
}
\dkz
By Lemma 4,  $f=f''f'$ and $g=g''g'$, with $f',f'',g',g''$
of the given form.
The terms $f'$ and $g'$ are completely determined by the codomain $B$
(the number of occurrences of each letter in $B$
determines $f'$ and $g'$) and therefore $f'$ and $g'$ are identical.
Suppose its type is $A \str A'$. Then the terms $f''$ and $g''$
are of the same type $A' \str B$.
Let $\alpha$ and $\beta$ be the {\bf SyMon} canonical transformations
corresponding to these terms.
Taking $\Gamma_{\alpha}$ and $\Gamma_{\beta}$
(the graphs of these transformations) as connections that connect an
occurrence of a letter
in $B$ with occurrence of the same letter  in $A'$ and since all
\mc-products in $f''$ and $g''$ are atomic diversified, they must connect
the first (from the left) occurrence of one letter in $B$ with the first
(again from the left) occurrence of the same letter in $A'$, the second with
the second etc. This means that $\Gamma_{\alpha}$ and $\Gamma_{\beta}$
are the same graphs and by coherence in symmetric monoidal categories,
$\alpha$ and $\beta$ are  the same canonical transformations in {\bf SyMon}.
Therefore, $f''$ and $g''$ are equal in {\bf SyMon}, and since all
symmetric monoidal equalities hold in {\bf Rel},
they are equal in {\bf Rel} too.
Hence, $f=g$ in {\bf Rel}.
\qed
Now Theorem 1  follows from lemmata 2 and 5.

\section{Coherence in Affine Categories}
As in the case of relevant categories, first we prove a lemma
about representation of {\bf Aff} morphism terms.
\\[0.3cm]
{\bf Lemma 6}\hspace{2em}
{\em
Every {\bf Aff} morphism-term is equal to a term of the form
$h_{2}h_{1}$, where $h_{1}$ is an ordered composition of atomic
\mk-products, and no \mk-product occurs in
$h_{2}$.
}
\dkz
As in the proof of Lemma 4, we transform the {\bf Aff}
morphism-term $h$ in several steps.
\\[0.3cm]
$1^{\circ}${\hspace{0.5cm}}
By the functoriality of multiplication, $h$ is equal to
a composition of product terms.
\\[0.3cm]
$2^{\circ}${\hspace{0.5cm}}
By the equalities (\mj\mk) and $\mk_{A{\cdot}B}=
\ms_{\mi}(\mj_{\mi}{\cdot}\mk_{B})(\mk_{A}{\cdot}\mj_{B})$,
which is derivable from (\mk) and (\mj\mk), 
this term is equal to a composition of
products with all \mk-products atomic.
\\[0.3cm]
$3^{\circ}${\hspace{0.5cm}}
By the naturality of \ms,\md,{\bf b},\mc-products and 
functoriality of multiplication, atomic \mk-products
permute to the right in order to obtain a term whose initial
part consists of all  \mk-products present in the term.
\\[0.3cm]
$4^{\circ}${\hspace{0.5cm}}
By the functoriality of multiplication, we can order this initial part
to obtain a term in the desired form, which is equal to $h$.
\qed
{\bf Lemma 7}\hspace{2em}
{\em
Let $f,g:A \str B$ be two morphism terms from {\bf Aff} with
$A$ diversified. Then $f=g$ in {\bf Aff}.
}
\dkz
By Lemma 6, there are terms $f',f'',g',g''$  such that
$f=f''f'$ and $g=g''g'$, where $f'$ and $g'$ are ordered compositions
of atomic \mk-products, and $f''$ and $g''$ are {\bf SyMon} terms.
The objects $A$ and $B$ (the letters occurring in $A$ not in $B$)
completely determine terms $f'$ and $g'$, hence they are identical
morphism terms, suppose of the type $A \str A'$.
Therefore
$f''$ and $g''$ are of the type $A' \str B$. By the assumption
concerning $A$ it
follows that $A'$ and $B$ are diversified.
Let $\alpha$ and $\beta$ be the canonical transformations in {\bf SyMon}
corresponding to $f''$ and $g''$ respectively. Taking $\Gamma_{\alpha}$ and
$\Gamma_{\beta}$ (theirs graphs) as connections between letter occurrences in
$B$ with letter occurrences in $A'$, since they connect a letter occurrence
with an occurrence of the same letter, and by the assumption about $A'$
and $B$, we must have that
$\Gamma_{\alpha}=\Gamma_{\beta}$. This implies, by the coherence
in symmetric monoidal categories, that $\alpha=\beta$, which has as 
consequence that $f''=g''$ in {\bf SyMon}, hence in {\bf Aff}.
We conclude that $f=f''f'=g''g'=g$.
\qed
From lemmata 2 and 7 it follows that
\\[0.3cm]
{\bf Theorem 2}\hspace{2em}
{\em
Every affine category is coherent.
}

\section{Coherence in Cartesian Categories}
The coherence in cartesian categories is not a new result. For the first time
it was mentioned in \cite{aaatc} and more recently in \cite{MI} and \cite{bpt}. Since we
would like to keep to the definition of coherence given above, we give 
another proof of this result here. 

One could expect that the proof of the coherence in cartesian
categories will follow the proofs given in the last two sections.
However, this method turns out to be too complicated and we use the standard
equational axiomatization of cartesian categories to avoid
this.
\\[0.3cm]
Denote by $\cal P$ the set of translations of morphism terms
from {\bf Cart} into the language of standard axiomatization
(see Section 1).
\\[0.3cm]
{\em Distributed terms} form the smallest class of morphism
terms from $\cal P$ that satisfies:
\begin{enumerate}
\item For all {\bf Cart} objects $A,B,C,D,E$, the term
$\mj_{A}$ as well as well-founded compositions of
$\mpm_{A,B},\mpm_{C,D}',\mk_{E}$ are in the class
and we call them {\em compat}.
\item If $f:C \str A$ and $g:C \str B$ are in the class
then ${\langle}f,g{\rangle}$ is in the class.
\end{enumerate}
The following  corresponds to the notion of the expanded normal form of
a natural deduction proof.
\\[0.2cm]
A distributed term is {\em atomic} if every compat
in this term has an atomic codomain (letter or \ri).
\\[0.3cm]
{\bf Lemma 8}\hspace{2em}
{\em
Every morphism-term from $\cal P$ is equal to an atomic distributed term.
}
\dkz
First, we show by induction on complexity of $f$ from $\cal P$
that it is equal to a distributed term.
\\[0.3cm]
$1^{\circ}${\hspace{0.5cm}}
If $f$ is $\mj_{A}$, $\mpm_{A,B}$,
$\mpm_{A,B}'$ or $\mk_{A}$, then it is compat and therefore distributed.
\\[0.3cm]
$2^{\circ}${\hspace{0.5cm}}
{\em a)} Suppose that $f$ is of the form ${\langle}g,h{\rangle}$.
By the induction hypothesis, $g$ and $h$ are equal to
distributed terms $g'$ and $h'$; hence $f$ is equal to the
distributed term ${\langle}g',h'{\rangle}$.
\\
{\em b)} Suppose that $f$ is of the form $hg$. Then by
the induction hypothesis, $g$ and $h$ are equal to distributed terms
$g_{1}$ and $h_{1}$. Suppose that every composition of lower complexity
than $h_{1}g_{1}$, of distributed terms 
is equal to a distributed term (if $g_{1}$ and $h_{1}$ are primitive,
then its composition is compat and therefore distributed). There are
three possibilities.
\\[0.2cm]
$i)${\hspace{0.5cm}}
If $g_{1}$ and $h_{1}$ are compat, then $h_{1}g_{1}$ is compat too,
hence $f$ is equal to the distributed term $h_{1}g_{1}$.
\\[0.2cm]
$ii)${\hspace{0.5cm}}
If $h_{1}$ is of the form ${\langle}j,l{\rangle}$,
for $j$ and $l$ distributed, then
$f={\langle}j,l{\rangle}g_{1}={\langle}jg_{1},lg_{1}{\rangle}$.
The terms $jg_{1}$ and $lg_{1}$ are of lower complexity than
$h_{1}g_{1}$, and
by assumption they are equal to some distributed
terms; hence $f$ is equal to a distributed term.
\\[0.2cm]
$iii)${\hspace{0.5cm}}
Suppose that $h_{1}$ is compat and $g_{1}$
is of the form ${\langle}j,l{\rangle}$.
If $h_{1}{\equiv}\mj$, then $f$ is equal to the distributed
term $g_{1}$. If $h_{1}{\equiv}h_{2}\mpm$, then
$f=h_{2}\mpm{\langle}j,l{\rangle}=h_{2}j$, where $h_{2}$ and $j$
are distributed and $h_{2}j$ is of
lower complexity than $h_{1}g_{1}$; therefore it is equal
to a distributed term. The case when
$h_{1}{\equiv}h_{2}\mpm'$ is analogous.
If $h_{1}{\equiv}h_{2}\mk$, then $f=h_{2}\mk{\langle}j,l{\rangle}=
h_{2}\mk$. The last term is compat, hence distributed.
\\[0.2cm]
This is the end of the induction. It follows that every $f$ from $\cal P$
is equal to a distributed term $f_{1}$. For every nonatomic
compat $h$ in $f_{1}$, using the equality
$h={\langle}\mpm_{A,B}h,\mpm_{A,B}'h{\rangle}$ for
$h:C \str A{\cdot}B$, we can find a distributed term equal to $h$,
such that every compat in it has a codomain of lower complexity than $h$.
Substituting this term for $h$ and repeating the procedure we obtain an
atomic distributed term equal to $f_{1}$ and therefore to $f$.
\qed
{\bf Lemma 9}\hspace{2em}
{\em
If $f,g:A \str B$ are two morphism-terms from $\cal P$ with
$A$ diversified, and $B$ an atom,
then $f=g$.
}
\dkz
If $B{\equiv}\ri$, then because it is terminal  in {\bf Cart},
we have that $f=g$. Suppose that $B{\equiv}q$ for $q$ a letter.
By the previous lemma $f$ and $g$ are equal to distributed terms $f_{1}$
and $g_{1}$, which are compat by the assumption that codomain is $q$.
Keeping in mind that there is no morphism in {\bf Cart} of the type
$A \str q$ such that $q$ doesn't occur in $A$, we prove the lemma
by induction on the complexity of the domain $A$.
\\[0.3cm]
$1^{\circ}${\hspace{0.5cm}}
If $A$ is an atom then $f_{1}{\equiv}g_{1}{\equiv}\mj_{q}$.
\\[0.3cm]
$2^{\circ}${\hspace{0.5cm}}
Suppose that $A{\equiv}A_{1}{\cdot}A_{2}$. Then, by the assumption,
$q$ occurs either in $A_{1}$ or in $A_{2}$. Suppose it occurs in
$A_{1}$. Then we must have that $f_{1}{\equiv}f_{2}\mpm$ and
$g_{1}{\equiv}g_{2}\mpm$ for some compat $f_{2},g_{2}:A_{1} \str q$.
By the induction hypothesis, $f_{2}=g_{2}$ holds, and therefore
$f=f_{1}=f_{2}\mpm=g_{2}\mpm=g_{1}=g$. We prove analogously the case when
$q$ occurs only in $A_{2}$.
\qed
{\bf Lemma 10}\hspace{2em}
{\em
Let $f,g:A \str B$ be two morphism-terms from $\cal P$
with $A$ diversified. Then $f=g$.
}
\dkz
By Lemma 8, $f$ and $g$ are equal to atomic distributed terms $f_{1}$
and $g_{1}$. The proof follows by induction on the complexity
of the codomain $B$.
If $B$ is an atom, then by the previous lemma $f=g$ holds.
If $B$ is not an atom, then neither $f_{1}$ nor $g_{1}$ are
compat, and therefore $f_{1}{\equiv}{\langle}i,j{\rangle}$ and
$g_{1}{\equiv}{\langle}l,h{\rangle}$, where $B=B_{1}{\cdot}B_{2}$.
The terms $i,l:A \str B_{1}$ and $j,h:A \str B_{2}$ are atomic
distributed and $B_{1}$, $B_{2}$ are of lower complexity than $B$.
Therefore, by the induction hypothesis, $i=l$ and $j=h$; hence $f=g$.
\qed
{\bf Corollary}\hspace{2em}
{\em
Let $f,g:A \str B$ be morphism terms from {\bf Cart} and let  $A$
be diversified. Then $f=g$ in {\bf Cart}.}
\dkz
Let {\boldmath $f$} and {\boldmath $g$} from $\cal P$ correspond
to $f$ and $g$ respectively. By Lemma 10,
{\boldmath $f$}$=${\boldmath $g$} and by the
extensional equivalence of these two axiomatizations we have that
$f=g$ in {\bf Cart}.
\qed
From this corollary and Lemma 2 it follows that
\\[0.3cm]
{\bf Theorem 3}\hspace{2em}
{\em
Every cartesian category is coherent.}

\section{Some consequences of the coherence}

\noindent
Usually, in the literature, coherence is not related to 
matters concerning natural transformations, but to conditions
that imply equality of morphisms of certain categories.
In the previous sections, we have lemmata 3, 5 and 7 as
examples of such an opinion. Now we prove some facts
that can be of practical interest for substructural categories,
especially for the free categories of each type.

In Section 4, a definition of canonical transformation that corresponds
to a morphism term $f$ of a free category $\cal C$ is given. From now on,
a graph of this transformation will be called a graph of $f$. By a 
straightforward induction on the complexity of $f$ we can prove the
following.
\\[0.3cm]
{\bf Lemma 11}\hspace{2em}
{\em If $f$ and $g$ are two equal morphism terms in $\cal C$, then 
the graph of $f$  is identical to the graph of $g$.}
\\[0.3cm]

Denote by ${\mbox{\bf Finord}}^{op}$ the dual of the category 
whose objects are finite
ordinals and whose arrows are mappings between them. The coherence of
substructural categories together with Lemma 11
is equivalent to the fact that there exist 
embeddings of  {\bf SyMon}, {\bf Rel}, {\bf Aff} and {\bf Cart} 
into ${\mbox{\bf Finord}}^{op}$ given by the ``graph''
functor $G$, such that for every $A \in {\cal O}$,
$G(A)$ is the number of occurrences of letters  in $A$ and 
$G(f)$ is the graph of $f$.

In the case of {\bf Cart}, this embedding is onto on morphisms,
and if we restrict ourserlves to
{\em one-one} functions in ${\mbox{\bf Finord}}^{op}$, then the embedding of
{\bf Aff} in this category is also onto on morphisms. Similarly we obtain
embeddings which are onto on morphisms
of {\bf Rel} and {\bf SyMon} in ${\mbox{\bf Finord}}^{op}$ with restrictions
to {\em onto} functions and {\em bijections} respectively.
This is an alternative characterization of the coherence of substructural
categories.

The results obtained in previous sections imply that
the categories {\bf Rel}, {\bf Aff} and {\bf Cart} are trivial in some
sense. However, they are not preorders (in a preorder, 
there is at most one arrow between two objects) as the case is with {\bf Mon},
though lemmata 5, 7 and 10 come close to preordering.
We can use the following consequence of coherence in these categories:
\\[0.2cm]

{\em Let $\cal C$ be one of the mentioned free categories and let
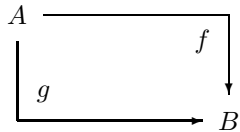
\begin{center}
\begin{picture}(120,60)
\put(20,10){\vector(1,0){70}}
\put(30,50){\line(1,0){70}}
\put(20,10){\line(0,1){30}}
\put(100,50){\vector(0,-1){30}}
\put(20,50){\makebox(0,0){$A$}}
\put(100,10){\makebox(0,0){$B$}}
\put(90,40){\makebox(0,0){$f$}}
\put(30,20){\makebox(0,0){$g$}}
\end{picture}
\end{center}
be one of its diagrams. It commutes iff the graphs of $f$ and $g$ 
are identical.}
\\[0.2cm]
Together with freedom of $\cal C$, this consequence can be of
practical use because it transforms computations in algebra of
morphism terms to the simple calculus of morphism graphs.
\\[0.3cm]
{\bf Example}{\hspace{2em}}
Suppose we want to simplify the term 
\[ ((\md_{A}(\mj_{A}{\cdot}\ms_{\mi})){\cdot}(\md_{B}\ms_{B{\cdot}\mi}))
\mc^{m}_{A,\mi,\mi{\cdot}\mi,B{\cdot}\mi}(\mj_{A{\cdot}\mi}\mc^{m}
_{\mi,B,\mi,\mi})(\md^{i}_{A}{\cdot}(\ms^{i}_{B}{\cdot}\md^{i}_{\mi}))
\]
representing a morphism of some
symmetric monoidal category $\cal S$.
The type of this term is
$A{\cdot}(B{\cdot}\ri){\str}A{\cdot}B$. 
Consider the {\bf SyMon} term
\[ ((\md_{p}(\mj_{p}{\cdot}\ms_{\mi})){\cdot}(\md_{q}\ms_{q{\cdot}\mi}))
\mc^{m}_{p,\mi,\mi{\cdot}\mi,q{\cdot}\mi}(\mj_{p{\cdot}\mi}\mc^{m}
_{\mi,q,\mi,\mi})(\md^{i}_{p}{\cdot}(\ms^{i}_{q}{\cdot}\md^{i}_{\mi})):
p{\cdot}(q{\cdot}\ri){\str}p{\cdot}q.
\]
Since its domain is diversified, it is enough to find a simple
{\bf SyMon} term of the same type. In this case, the term
$\mj_{p}{\cdot}\md_{q}$ is imposed. So, these two terms are equal
in {\bf SyMon}, and by the freedom of this category, 
the initial term is equal to
$\mj_{A}{\cdot}\md_{B}$ in $\cal S$.
\\[0.3cm]
{\bf Example}{\hspace{2em}}
Prove the equality of the following {\bf Cart} terms
\[(\mj_{p}{\cdot}\mc_{p,q{\cdot}p})\bl_{p,p,q{\cdot}p}(\mj_{p{\cdot}p}
{\cdot}\mc_{p,q})(\mw_{p}{\cdot}\mj_{p{\cdot}q}),\]
and
\[(\mj_{p}{\cdot}(((\md_{q{\cdot}p}\br_{q,p,\mi}){\cdot}\mj_{p})\br_{q,
p{\cdot}\mi,p}))\bl_{p,q,(p{\cdot}\mi){\cdot}p}((\mj_{p}{\cdot}\ms_{q})
{\cdot}((\mj_{p}{\cdot}\mk_{q}){\cdot}\mj_{p}))((\mj_{p}{\cdot}(
\mk_{p}{\cdot}\mj_{q})){\cdot}
\mc_{p,p{\cdot}q})\mw_{p{\cdot}(p{\cdot}q)}.\]
Their graphs are identical, which can be checked directly from
the constructions of these terms, following the linkages between
the  occurrences of letters in domains and codomains of its primitive
components.
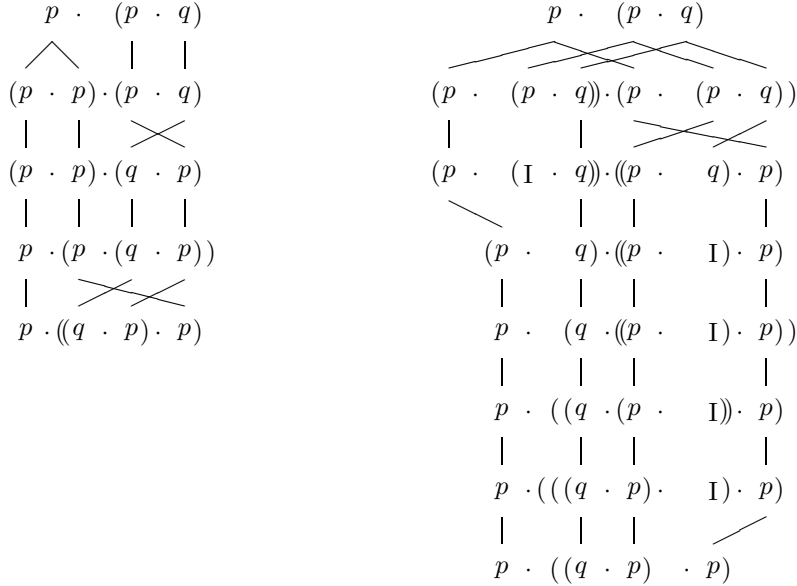
\begin{center}
\rm
\begin{picture}(300,250)
\put(30,230){\circle*{1}}
\put(60,230){\circle*{1}}
\put(20,200){\circle*{1}}
\put(40,200){\circle*{1}}
\put(60,200){\circle*{1}}
\put(20,170){\circle*{1}}
\put(40,170){\circle*{1}}
\put(60,170){\circle*{1}}
\put(20,140){\circle*{1}}
\put(40,140){\circle*{1}}
\put(60,140){\circle*{1}}
\put(18,110){\circle*{1}}
\put(40,110){\circle*{1}}
\put(60,110){\circle*{1}}
\put(20,230){\makebox(0,0){$p$}}
\put(50,230){\makebox(0,0){$p$}}
\put(70,230){\makebox(0,0){$q$}}
\put(10,200){\makebox(0,0){$p$}}
\put(30,200){\makebox(0,0){$p$}}
\put(50,200){\makebox(0,0){$p$}}
\put(70,200){\makebox(0,0){$q$}}
\put(10,170){\makebox(0,0){$p$}}
\put(30,170){\makebox(0,0){$p$}}
\put(50,170){\makebox(0,0){$q$}}
\put(70,170){\makebox(0,0){$p$}}
\put(10,140){\makebox(0,0){$p$}}
\put(30,140){\makebox(0,0){$p$}}
\put(50,140){\makebox(0,0){$q$}}
\put(70,140){\makebox(0,0){$p$}}
\put(10,110){\makebox(0,0){$p$}}
\put(30,110){\makebox(0,0){$q$}}
\put(50,110){\makebox(0,0){$p$}}
\put(70,110){\makebox(0,0){$p$}}
\put(20,220){\line(-1,-1){10}}
\put(20,220){\line(1,-1){10}}
\put(50,220){\line(0,-1){10}}
\put(70,220){\line(0,-1){10}}
\put(10,190){\line(0,-1){10}}
\put(30,190){\line(0,-1){10}}
\put(50,190){\line(2,-1){20}}
\put(70,190){\line(-2,-1){20}}
\put(10,160){\line(0,-1){10}}
\put(30,160){\line(0,-1){10}}
\put(50,160){\line(0,-1){10}}
\put(70,160){\line(0,-1){10}}
\put(10,130){\line(0,-1){10}}
\put(30,130){\line(4,-1){40}}
\put(50,130){\line(-2,-1){20}}
\put(70,130){\line(-2,-1){20}}
\put(45,230){\makebox(0,0){$($}}
\put(75,230){\makebox(0,0){$)$}}
\put(5,200){\makebox(0,0){$($}}
\put(35,200){\makebox(0,0){$)$}}
\put(45,200){\makebox(0,0){$($}}
\put(75,200){\makebox(0,0){$)$}}
\put(5,170){\makebox(0,0){$($}}
\put(35,170){\makebox(0,0){$)$}}
\put(45,170){\makebox(0,0){$($}}
\put(75,170){\makebox(0,0){$)$}}
\put(25,140){\makebox(0,0){$($}}
\put(45,140){\makebox(0,0){$($}}
\put(75,140){\makebox(0,0){$)$}}
\put(80,140){\makebox(0,0){$)$}}
\put(23,110){\makebox(0,0){$($}}
\put(25,110){\makebox(0,0){$($}}
\put(55,110){\makebox(0,0){$)$}}
\put(75,110){\makebox(0,0){$)$}}
\put(220,230){\circle*{1}}
\put(250,230){\circle*{1}}
\put(180,200){\circle*{1}}
\put(210,200){\circle*{1}}
\put(230,200){\circle*{1}}
\put(250,200){\circle*{1}}
\put(280,200){\circle*{1}}
\put(180,170){\circle*{1}}
\put(210,170){\circle*{1}}
\put(230,170){\circle*{1}}
\put(250,170){\circle*{1}}
\put(280,170){\circle*{1}}
\put(200,140){\circle*{1}}
\put(230,140){\circle*{1}}
\put(250,140){\circle*{1}}
\put(280,140){\circle*{1}}
\put(200,110){\circle*{1}}
\put(230,110){\circle*{1}}
\put(250,110){\circle*{1}}
\put(280,110){\circle*{1}}
\put(200,80){\circle*{1}}
\put(230,80){\circle*{1}}
\put(250,80){\circle*{1}}
\put(280,80){\circle*{1}}
\put(200,50){\circle*{1}}
\put(230,50){\circle*{1}}
\put(250,50){\circle*{1}}
\put(280,50){\circle*{1}}
\put(200,20){\circle*{1}}
\put(230,20){\circle*{1}}
\put(260,20){\circle*{1}}
\put(210,230){\makebox(0,0){$p$}}
\put(240,230){\makebox(0,0){$p$}}
\put(260,230){\makebox(0,0){$q$}}
\put(170,200){\makebox(0,0){$p$}}
\put(200,200){\makebox(0,0){$p$}}
\put(220,200){\makebox(0,0){$q$}}
\put(240,200){\makebox(0,0){$p$}}
\put(270,200){\makebox(0,0){$p$}}
\put(290,200){\makebox(0,0){$q$}}
\put(170,170){\makebox(0,0){$p$}}
\put(200,170){\makebox(0,0){I}}
\put(220,170){\makebox(0,0){$q$}}
\put(240,170){\makebox(0,0){$p$}}
\put(270,170){\makebox(0,0){$q$}}
\put(290,170){\makebox(0,0){$p$}}
\put(190,140){\makebox(0,0){$p$}}
\put(220,140){\makebox(0,0){$q$}}
\put(240,140){\makebox(0,0){$p$}}
\put(270,140){\makebox(0,0){I}}
\put(290,140){\makebox(0,0){$p$}}
\put(190,110){\makebox(0,0){$p$}}
\put(220,110){\makebox(0,0){$q$}}
\put(240,110){\makebox(0,0){$p$}}
\put(270,110){\makebox(0,0){I}}
\put(290,110){\makebox(0,0){$p$}}
\put(190,80){\makebox(0,0){$p$}}
\put(220,80){\makebox(0,0){$q$}}
\put(240,80){\makebox(0,0){$p$}}
\put(270,80){\makebox(0,0){I}}
\put(290,80){\makebox(0,0){$p$}}
\put(190,50){\makebox(0,0){$p$}}
\put(220,50){\makebox(0,0){$q$}}
\put(240,50){\makebox(0,0){$p$}}
\put(270,50){\makebox(0,0){I}}
\put(290,50){\makebox(0,0){$p$}}
\put(190,20){\makebox(0,0){$p$}}
\put(220,20){\makebox(0,0){$q$}}
\put(240,20){\makebox(0,0){$p$}}
\put(270,20){\makebox(0,0){$p$}}
\put(210,220){\line(-4,-1){40}}
\put(240,220){\line(-4,-1){40}}
\put(260,220){\line(-4,-1){40}}
\put(210,220){\line(3,-1){30}}
\put(240,220){\line(3,-1){30}}
\put(260,220){\line(3,-1){30}}
\put(170,190){\line(0,-1){10}}
\put(220,190){\line(0,-1){10}}
\put(240,190){\line(5,-1){50}}
\put(270,190){\line(-3,-1){30}}
\put(290,190){\line(-2,-1){20}}
\put(170,160){\line(2,-1){20}}
\put(220,160){\line(0,-1){10}}
\put(240,160){\line(0,-1){10}}
\put(290,160){\line(0,-1){10}}
\put(190,130){\line(0,-1){10}}
\put(220,130){\line(0,-1){10}}
\put(240,130){\line(0,-1){10}}
\put(290,130){\line(0,-1){10}}
\put(190,100){\line(0,-1){10}}
\put(220,100){\line(0,-1){10}}
\put(240,100){\line(0,-1){10}}
\put(290,100){\line(0,-1){10}}
\put(190,70){\line(0,-1){10}}
\put(220,70){\line(0,-1){10}}
\put(240,70){\line(0,-1){10}}
\put(290,70){\line(0,-1){10}}
\put(190,40){\line(0,-1){10}}
\put(220,40){\line(0,-1){10}}
\put(240,40){\line(0,-1){10}}
\put(290,40){\line(-2,-1){20}}
\put(235,230){\makebox(0,0){$($}}
\put(265,230){\makebox(0,0){$)$}}
\put(165,200){\makebox(0,0){$($}}
\put(195,200){\makebox(0,0){$($}}
\put(224,200){\makebox(0,0){$)$}}
\put(226,200){\makebox(0,0){$)$}}
\put(235,200){\makebox(0,0){$($}}
\put(265,200){\makebox(0,0){$($}}
\put(295,200){\makebox(0,0){$)$}}
\put(300,200){\makebox(0,0){$)$}}
\put(165,170){\makebox(0,0){$($}}
\put(195,170){\makebox(0,0){$($}}
\put(224,170){\makebox(0,0){$)$}}
\put(226,170){\makebox(0,0){$)$}}
\put(234,170){\makebox(0,0){$($}}
\put(236,170){\makebox(0,0){$($}}
\put(275,170){\makebox(0,0){$)$}}
\put(295,170){\makebox(0,0){$)$}}
\put(185,140){\makebox(0,0){$($}}
\put(225,140){\makebox(0,0){$)$}}
\put(234,140){\makebox(0,0){$($}}
\put(236,140){\makebox(0,0){$($}}
\put(275,140){\makebox(0,0){$)$}}
\put(295,140){\makebox(0,0){$)$}}
\put(234,110){\makebox(0,0){$($}}
\put(215,110){\makebox(0,0){$($}}
\put(236,110){\makebox(0,0){$($}}
\put(300,110){\makebox(0,0){$)$}}
\put(275,110){\makebox(0,0){$)$}}
\put(295,110){\makebox(0,0){$)$}}
\put(210,80){\makebox(0,0){$($}}
\put(215,80){\makebox(0,0){$($}}
\put(235,80){\makebox(0,0){$($}}
\put(276,80){\makebox(0,0){$)$}}
\put(274,80){\makebox(0,0){$)$}}
\put(295,80){\makebox(0,0){$)$}}
\put(210,50){\makebox(0,0){$($}}
\put(215,50){\makebox(0,0){$($}}
\put(205,50){\makebox(0,0){$($}}
\put(245,50){\makebox(0,0){$)$}}
\put(275,50){\makebox(0,0){$)$}}
\put(295,50){\makebox(0,0){$)$}}
\put(210,20){\makebox(0,0){$($}}
\put(215,20){\makebox(0,0){$($}}
\put(245,20){\makebox(0,0){$)$}}
\put(275,20){\makebox(0,0){$)$}}
\end{picture}
\end{center}

Since these terms are of the same type, they are equal in 
{\bf Cart}.

We conclude this section with a discussion about 
the hierarchy (in the sense of embeddability)
in the set of free categories mentioned above.
All these categories have the same set $\cal O$ as the set of objects, 
and their morphism
terms satisfy the following inclusions
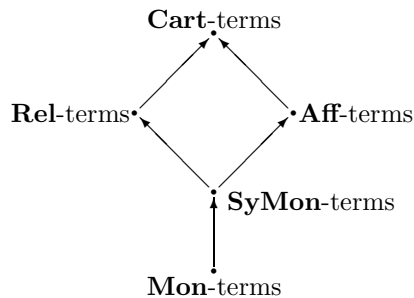
\begin{center}
\begin{picture}(160,130)
\put(80,20){\circle*{2}}
\put(80,18){\makebox(0,0)[t]{{\bf Mon}-terms}}
\put(80,50){\circle*{2}}
\put(85,50){\makebox(0,0)[tl]{{\bf SyMon}-terms}}
\put(50,80){\circle*{2}}
\put(48,80){\makebox(0,0)[r]{{\bf Rel}-terms}}
\put(110,80){\circle*{2}}
\put(112,80){\makebox(0,0)[l]{{\bf Aff}-terms}}
\put(80,110){\circle*{2}}
\put(80,112){\makebox(0,0)[b]{{\bf Cart}-terms}}
\put(80,22){\vector(0,1){26}}
\put(82,52){\vector(1,1){26}}
\put(78,52){\vector(-1,1){26}}
\put(108,82){\vector(-1,1){26}}
\put(52,82){\vector(1,1){26}}
\end{picture}
\end{center}

where arrows stay for inclusions. The fact that
equalities between morphism terms that hold in the ``lower'' category
are also true in the ``higher'' one was used in this paper several times.
To show that a ``lower'' category
is a subcategory of those above, in the diagram,
we have to show the following
\\[0.3cm]
{\bf Lemma 12}\hspace{2em}
{\em
Let $\cal C$ and $\cal D$ be categories from the above diagram
such that $\cal C$ is below $\cal D$.
If for $\cal C$ morphism terms $f$ and $g$, the equality $f=g$ holds
in $\cal D$, then they are equal in $\cal C$, too.
}
\dkz
Let $\Phi$ and $\Psi$ be the graphs of $f$ and $g$ respectively.
Since $f=g$ in $\cal D$, by Lemma 11, the
graphs $\Phi$ and $\Psi$
are identical. By the coherence of $\cal C$ (there is an unique
morphism of a certain type in a canonical transformation
of the free category $\cal C$) it follows that $f=g$ in $\cal C$.
\qed
So, by the embedding $E$ such that $E(A)=A$ for every $A \in {\cal O}$, and
$E([f]_{\cal C})=[f]_{\cal D}$ where $[f]_{\cal C}$ is the equivalence
class of the morphism term $f$ in $\cal C$,  we have the following
\\[0.3cm]
{\bf Theorem}\hspace{2em}
{\em
The category {\bf Mon} is a subcategory of {\bf SyMon}, {\bf Rel}, {\bf Aff}
and {\bf Cart}. The category {\bf SyMon} is a subcategory of {\bf Rel}, {\bf Aff}
and {\bf Cart}. The categories {\bf Rel} and {\bf Aff} are subcategories of
{\bf Cart}.}
\\[0.3cm]
This theorem looks almost trivial, but it is of the same strength as
the coherence theorems of substructural categories.
An independent proof of this theorem immediately
delivers MacLane's coherence for monoidal and symmetric monoidal
categories and the above coherence theorems for relevant and affine
categories, from the simplest case of cartesian coherence.

Some other coherence applications in this spirit are given in \cite{tmocc},
and the main result of \cite{ciasm} can be proven
more easily by using this apparatus.
\\[0.7cm]
{\it Acknowledgements.} I would like to thank Professor Kosta Do\v sen
for the inspiration and lot of helpful comments on this paper.
This work was supported by Grant 0401A of the Science Fund
of Serbia.

\vskip 0.25in
\noindent Zoran Petri\' c\\
University of Belgrade\\
Department of Mining and Geology\\
Dju\v sina 7\\
11000 Belgrade, Yugoslavia\\
e-mail: zpetric@rgf.rgf.bg.ac.yu
\pagebreak

\end{document}